\documentclass[final,3p,times]{elsarticle}
\usepackage{hyperref}

\journal{}

\usepackage{amssymb}
\usepackage{amsthm}
\usepackage{amsmath}
\usepackage{graphicx}
\usepackage{float}
\usepackage{subfigure}
\usepackage{caption}
\newcommand{\eqnref}[1]{(\ref{#1})}

\begin{document}

\begin{frontmatter}

\title{An Entropy Stable Central Solver for Euler Equations}




\author{N.H. Maruthi\fnref{fn1}}
\ead{maruthi@aero.iisc.ernet.in, maruthinh@gmail.com}

\author{S.V. Raghurama Rao\corref{cor1}\fnref{fn2}}
\ead{raghu@aero.iisc.ernet.in, svraghuramarao@gmail.com}

\cortext[cor1]{Corresponding author}
\fntext[fn1]{Graduate student}
\fntext[fn2]{Associate professor}
\address{Department of Aerospace Engineering, Indian Institute of Science, Bengaluru, Karnataka, India-560012}

\begin{abstract}  
     An exact discontinuity capturing central solver developed recently, named MOVERS (Method of Optimal Viscosity for Enhanced Resolution of Shocks) \cite{jaisankar2009central}, is analyzed and improved further to make it entropy stable.  MOVERS, which is designed to capture steady shocks and contact discontinuities exactly by enforcing the Rankine-Hugoniot jump condition directly in the discretization process, is a low diffusive algorithm in a simple central discretization framework, free of complicated Riemann solvers and flux splittings.  However, this algorithm needs an entropy fix to avoid nonsmoothness in the expansion regions.  The entropy conservation equation is used as a guideline to introduce an optimal numerical diffusion in the smooth regions and a limiter based switchover is introduced for numerical diffusion based on jump conditions at the large gradients.  The resulting new scheme is entropy stable, accurate and captures steady discontinuities exactly while avoiding an entropy fix.
\end{abstract}


\end{frontmatter}


\section{Introduction}
\label{Introduction} 
     Of all the numerical methods developed to solve the hyperbolic systems of conservation laws of fluid dynamics numerically, approximate Riemann solvers have been the most popular.  However, these approximate Riemann solvers are discovered to introduce shock instabilities like carbuncle shocks, odd-even decoupling, kinked Mach stems in the solutions \cite{quirk1994a} and they typically require addition of numerical diffusion in the form of entropy fixes to make them entropy stable.  It is often nontrivial to extend these approximate Riemann solvers to hyperbolic systems other than the Euler equations of gas dynamics and further they are  heavily dependent on the eigen-structure of the hyperbolic systems.  This dependency is due to the enforcement of strict upwinding procedure which usually results in sophisticated procedures of Riemann solvers or complicated flux splittings and are dependent on the projection of the solution on to the space of eigenvectors and characteistic decomposition.  In this context, in recent research, attention is focused on developing simpler central solvers which are competitive to upwind schemes in terms of accuracy and efficiency of capturing essential flow features like shock waves, contact discontinuities and expansion waves \cite{Nessyahu_Tadmor, Kurganov_Tadmor, jaisankar2007diffusion, jaisankar2009central, central_station}.  

     Of the new central schemes, the {\em Method of Optimal Viscosity for Enhanced Resolution of Shocks (MOVERS)} developed by Jaisankar and Raghurama Rao \cite{jaisankar2009central} is capable of capturing steady discontinuities like shock waves and contact discontinuities aligned with grid lines exactly without any numerical diffusion and is otherwise a low numerical diffusion algorithm.  The  exact discontinuity capturing, a feature shared by very few of the upwind schemes which made them popular, is a result of enforcing the Rankine-Hugoniot jump condition directly in the discretization process in MOVERS.  With this feature the central discretization methods can compete in accuracy with the best of the Riemann solvers, in spite of being much simpler.  However, the two versions of MOVERS presented by Jaisankar and Raghurama Rao \cite{jaisankar2009central}, namely MOVERS-n (an n-wave model for the coefficent of numerical diffusion) and MOVERS-1 (a 1-wave model for the coefficient of numerical diffusion) require an entropy fix to avoid non-smoothness in the expansion regions.  To overcome the problem of an entropy fix, another version, MOVERS-L, is presented in that work based on the idea of shifting to base-line Rusanov solver \cite{rusanov1962calculation} or Local Lax-Friedrichs (LLF) method \cite{Chi-Wang_Shu_1998, leveque2002finite} in smooth regions and to MOVERS-n or MOVERS-1 near large gradients based on a limiter function.  It is worth noting that while MOVERS-L overcomes the problem of using an entropy fix, due to the base-line solver being LLF method, a high dose of numerical diffusion will be present in the smooth regions while still capturing the discontinuities very accurately.  In this work, we propose an alternative method to reduce the numerical diffusion in the smooth regions in MOVERS-L, based on appropriate guidance from the discretization of the entropy conservation equation.  For a detailed mathematical analysis of entropy stability the reader is referred to the research work of Tadmor \cite{tadmor2003entropy}.   In this work, however, we present a simpler way of introducing entropy stability based on an LLF type discretization of the entropy conservation equation and using it for guidance while fixing the numerical diffusion in smooth regions for MOVERS-L, thereby introducing a new version of the algorithm which we name as MOVERS-LE.  
			  
\section{Numerical Method}
     Before presenting the new entropy stable version of the algorithm, MOVERS-LE, its base-line solvers, MOVERS and MOVERS-L, are first introduced briefly in the following.  

\subsection{MOVERS for 1-D Euler equations}
Consider the 1-D Euler equations given by  
\begin{equation}\label{eq1}
 \frac{\partial{U}}{\partial{t}} + \frac{\partial{F}}{\partial{x}} = 0 
\end{equation}
where the conserved variable vector U and the inviscid flux vector F are defined as

\begin{equation*}
 U = \begin{bmatrix}
       \rho        \\[0.3em]
       \rho u \\[0.3em]
      \rho E
     \end{bmatrix}
\ \  \textrm{and} \ \
 F = \begin{bmatrix}
       \rho u       \\[0.3em]
        p + \rho u^2  \\[0.3em]
       pu + \rho u E 
     \end{bmatrix}
\end{equation*}  
Here, $\rho$ is the density, u is the velocity, E is the total energy, given by $E = \frac{p}{\rho(\gamma-1)} + \frac{u^2}{2}$ and p is the pressure. The flux Jacobian matrix  defined by $A = \frac{\partial{F}}{\partial{U}}$ has three real and distinct eigenvalues, $\lambda_1=u-a, \lambda_2 = u$ and$ \ \lambda_3 = u + a$. Applying a semi-discrete cell-centered finite volume method to this system of conservation laws gives the update formula \\
\begin{equation}\label{eq2}
   \frac{dU_j}{dt} =  -  \frac{1}{\Delta{x}}[ F^{n}_{j+\frac{1}{2}}-  F^{n}_{j-\frac{1}{2}}] 
\end{equation} 
where the superscript n indicates the time level and the subscript $j$ represents the grid points at the cell centres, with $j\pm\frac{1}{2}$ representing the cell interfaces. The cell interface flux is written as an average flux plus a diffusive flux given by
\begin{equation}\label{eq3}
 F^{n}_{j+\frac{1}{2}}= \frac{1}{2}(F_j + F_{j+1})- \frac{1}{2}\alpha_{j+\frac{1}{2}}(U_{j+1} - U_j) 
\end{equation}
where $\alpha_{i+\frac{1}{2}}$ represents the coefficient of numerical diffusion. 
For a stable scheme the general interface flux is given by  
\begin{equation} \label{eq4}
F_{I} = \frac{1}{2}(F_{L} + F_{R}) - \frac{1}{2} \alpha_I \Delta U
\end{equation}
where $L$ and $R$ represent the left and right states of the interface, the first term on the right hand side (average flux) represents the central discretization of the flux terms and the second term is the diffusive flux, with $\alpha_I$ being the coefficient of numerical diffusion at the interface. The intention here is to obtain $\alpha$ using R-H condition to resolve steady discontinuities exactly. In this setting the relevant physical speed is the speed of the discontinuity.  This speed, $s$,  is given by Rankine-Hugoniot jump condition. Here R-H condition is enforced at the cell interfaces to get the relevant expression for coefficient of numerical diffusion.  For a discontinuous wave with a left state and right state given by $U_{L}$ and $U_{R}$, the Rankine-Hugoniot jump condition is given by 

\begin{equation}\label{eq5}
\begin{split}
\Delta F & = s \Delta U \\
~~\mbox{where, }~~
\Delta F & = F_{R} - F_{L} \\
\Delta U & = U_{R} - U_{L}
\end{split}
\end{equation}

\begin{figure}
\begin{center}
\includegraphics[trim=30.0 0.0 50.0 10.0, clip, width=0.6\textwidth]{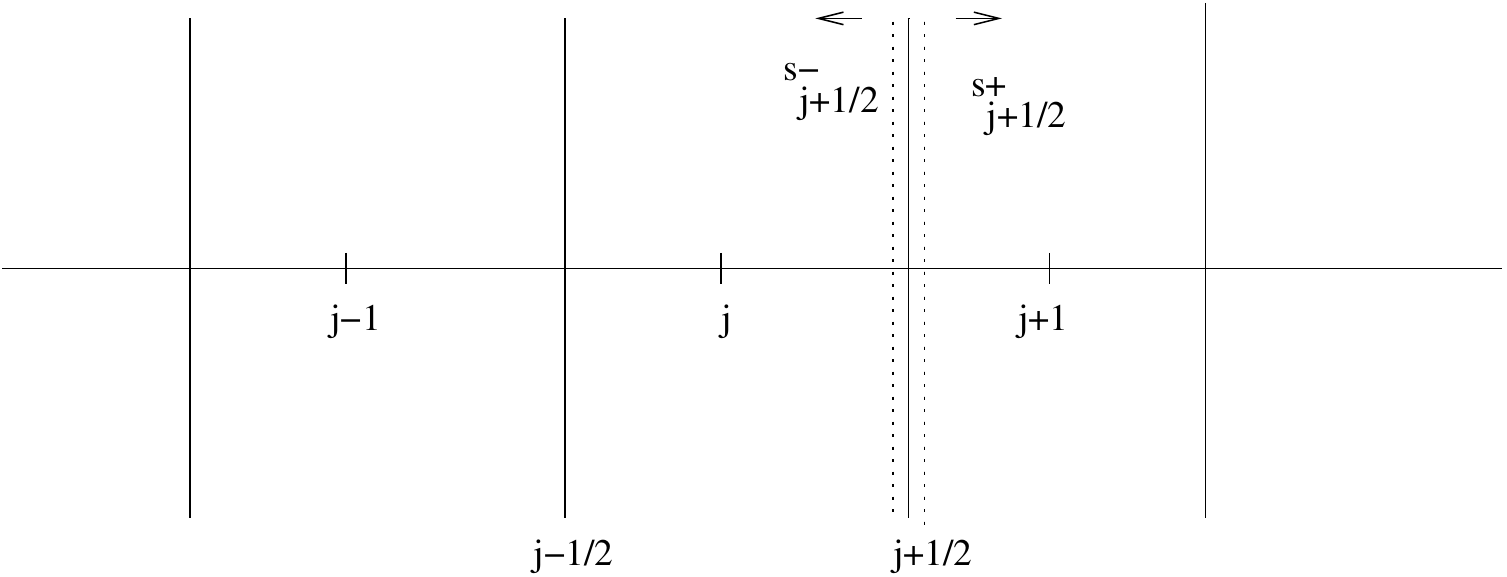}
\caption{Moving Shock }
\label{shock_moving}
\end{center}
\end{figure}

In the finite volume method, the variation of the conserved variable is assumed as a  piecewise polynomial function within each finite volume.  Hence the discontinuity in the conserved variable can occur at a cell interface.  Let us apply the above R-H condition at a cell interface, $j+\frac{1}{2}$. 
\begin{equation}\label{eq6}
F_{j+1} - F_{j} = s_{j+\frac{1}{2}} \left( U_{j+1} - U_{j} \right)
\end{equation}
To obtain the interface flux, we split the above R-H condition into two parts.  Let us consider a discontinuity between $j$ and  $j+1$ moving with a speed $s_{j+\frac{1}{2}}$ as shown in Fig~(\ref{shock_moving}).
This discontinuity can move to either left or right. The jump condition does not indicate the direction of the movement of the discontinuity. To distinguish the positive and negative speeds, we can split the speed of the discontinuity into two parts, positive and negative. By this, left and right moving discontinuities are separated.  We can now split the R-H condition also into two parts as 
\begin{equation} \label{RH1}
F_{j+1} - F_{j+\frac{1}{2}} = s^{+}_{j+\frac{1}{2}} \left( U_{j+1} - U_{j} \right) 
\end{equation} 
and 
\begin{equation} \label{RH2}
F_{j+\frac{1}{2}} - F_{j} = s^{-}_{j+\frac{1}{2}} \left( U_{j+1} - U_{j} \right) 
\end{equation} 
Note that the right hand side is not changed as in the above equation of three variables, namely $s$, $U$ and $F$, only two variables can be varied.  Subtracting (\ref{RH2}) from (\ref{RH1}), we obtain 
\begin{equation} 
F_{j+1} - 2 F_{j+\frac{1}{2}} - F_{j} = \left( s^{+}_{j+\frac{1}{2}} - s^{-}_{j+\frac{1}{2}} \right)
\left( U_{j+1} - U_{j} \right) 
\end{equation} 
or 
\begin{equation} \label{eq7} 
F_{j+\frac{1}{2}} = \frac{1}{2} \left[ F_{j} + F_{j+1} \right] 
- \left( \frac{s_{j+\frac{1}{2}}+|s_{j+\frac{1}{2}}|}{2} 
- \frac{s_{j+\frac{1}{2}}-|s_{j+\frac{1}{2}}|}{2} \right)  
\left[ U_{j+1} - U_{j} \right]   
\end{equation} 
Now we get the interface flux in-terms of the speed of the discontinuity which is obtained using R-H condition as
\begin{equation} \label{eq8}
F_{j+\frac{1}{2}} = \frac{1}{2}\left( F_{j} + F_{j+1} \right) - \frac{1}{2} \mid s_{j+\frac{1}{2}} \mid(U_{j+1}-U_{j})
\end{equation}
Comparing both the interface fluxes (\ref{eq3}) and (\ref{eq8}), we obtain 
\begin{equation}\label{eq9}
\alpha_{j+\frac{1}{2}} = \mid {s_{j+\frac{1}{2}}} \mid
\end{equation}
Thus the coefficient of numerical diffusion is equal to the absolute value of the speed of the discontinuity.  We can evaluate the value of $s_{j+\frac{1}{2}}$ numerically which will be described in the following.  

Note that R-H condition is used to obtain coefficient of numerical diffusion in MOVERS.  Thus, this scheme captures steady discontinuities (shock waves, contact discontinuities, slip surfaces) exactly. To ensure that the speed of the discontinuity falls within the range of the eigen-spectrum, wave-speed correction is devised, which is discussed in the paper~\cite{jaisankar2009central}.

R-H condition for a system of hyperbolic conservation laws is given by
\begin{equation}\label{eq11}
\Delta F = s \Delta U 
\end{equation} 
where 
\begin{equation} 
\Delta F = F_{R} - F_{L} \ \textrm{and} \ 
\Delta U = U_{R} - U_{L}
\end{equation} 
Here,for 1-D Euler equations, $\Delta F$  and $\Delta U$ are $3 \times 1$ matrices, hence we take $s_{j+\frac{1}{2}}$ to be a $3 \times 3$ matrix.  We choose $s_{j+\frac{1}{2}}$ to be a simple diagonal matrix as follows.

\begin{equation}\label{eq12}
\left(
  \begin{array}{ccc}
       \Delta F_1        \\[0.3em]
       \Delta F_2 \\[0.3em]
      \Delta F_3
     \end{array}
\right)=
\left(
\begin{array}{ccc}
 s_{1,j+\frac{1}{2}} & 0 & 0 \\
 0 & s_{2,j+\frac{1}{2}} & 0 \\
 0 & 0 & s_{3,j+\frac{1}{2}} \\
\end{array}
\right)
\left(
  \begin{array}{ccc}
       \Delta U_1        \\[0.3em]
       \Delta U_2 \\[0.3em]
      \Delta U_3
     \end{array}
\right)
\end{equation}
Thus we obtain 
\begin{equation}\label{eq13}
\left(
  \begin{array}{ccc}
       s_{1,j+\frac{1}{2}}        \\[0.3em]
       s_{2,j+\frac{1}{2}} \\[0.3em]
      s_{3,j+\frac{1}{2}}
     \end{array}
\right)=
\left(
  \begin{array}{ccc}
       \frac{\Delta F1}{\Delta U1}        \\[0.3em]
       \frac{\Delta F2}{\Delta U2} \\[0.3em]
      \frac{\Delta F3}{\Delta U3}
     \end{array}
\right)
\end{equation}
Note that the coefficient of numerical diffusion $\alpha_{j+\frac{1}{2}}$ is to be fixed as the absolute value of the speed of the discontinuity as we derived and we thus now have three speeds of discontinuities, each for one of the conservation laws.  Since we have $n$ number of wave speeds representing $n$ conservation laws as the coefficients of numerical diffusion, this algorithm is termed as MOVERS-n. This can be easily extended to any hyperbolic system of conservation laws. In smooth regions of the flow, $\Delta U \approx 0$ and the evaluated diffusion values can be unphysical. To avoid this numerical overflow, the wave speed correction strategy is devised. The wave-speed correction is done such that the  coefficient of numerical diffusion is limited within the range of the eigen-spectrum of the flux Jacobian matrix. In this framework the wave-speed correction is slightly modified. When the left flux is equal to right flux and there is jump in conserved variables, then this case refers to steady discontinuity, for which the speed of the discontinuity is zero and hence the steady shock waves and contact discontinuities are captured exactly without any numerical diffusion. Considering all the possibilities discussed, the implementation of coefficient of numerical diffusion is discussed below.  

\begin{equation}\label{rhc}
 s_{j+\frac{1}{2}} = \begin{cases}
  \left|\frac{F_{j+1} - F_{j}} {U_{j+1} - U_{j}}\right|
  ~~\mbox{if}~~\ \ U_{j} \neq U_{j+1} ~~\mbox{and}~~\ F_{j} \neq F_{j+1} \\ 
   0 
  ~~\mbox{if}~~\  U_{j} \neq U_{j+1}  ~~\mbox{and}~~\ \left|F_{j+1} - F_{j}\right| < \epsilon \\
	  \lambda_{min} 
  ~~\mbox{if}~~\  \left|U_{j+1} - U_{j}\right| < \epsilon \end{cases}
\end{equation}
where $\epsilon$ is a small parameter.  Wave speed correction is enforced as follows.  
\begin{eqnarray} \label{wavecor}
  s_{j+\frac{1}{2}}= \lambda_{max}  ~~\mbox{if}~~\ s_{j+\frac{1}{2}} > \lambda_{max} \\
  s_{j+\frac{1}{2}}= \lambda_{min}  ~~\mbox{if}~~\ s_{j+\frac{1}{2}} < \lambda_{min}
 \end{eqnarray}
Here, $\lambda_{max}$ and $\lambda_{min}$ are maximum and minimum eigenvalues computed across the interface and are given as
 \begin{eqnarray} \label{eigval}
 \lambda_{max} =  max\{max(|u+a|, |u|, |u-a|)_L, max(|u+a|, |u|, |u-a|)_R \} \\
\lambda_{min} =  max\{min(|u+a|, |u|, |u-a|)_L, min(|u+a|, |u|, |u-a|)_R \} 
\label{mineig}
\end{eqnarray}

\begin{figure}[h!]
 \subfigure{\includegraphics[trim=75.0 125.0 100.0 131.0, clip, width=0.5\textwidth]{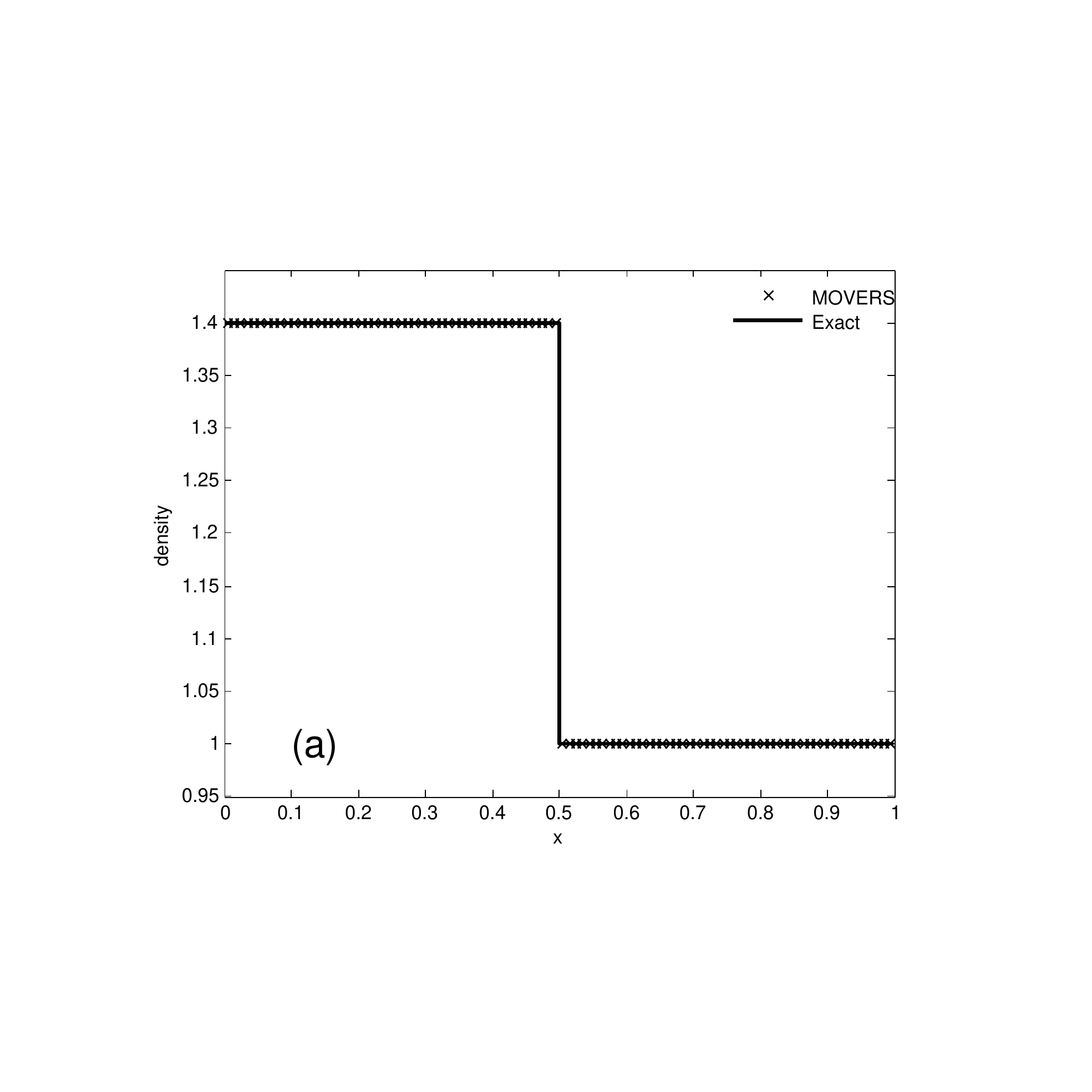}}
 \subfigure{\includegraphics[trim=75.0 125.0 100.0 131.0, clip, width=0.5\textwidth]{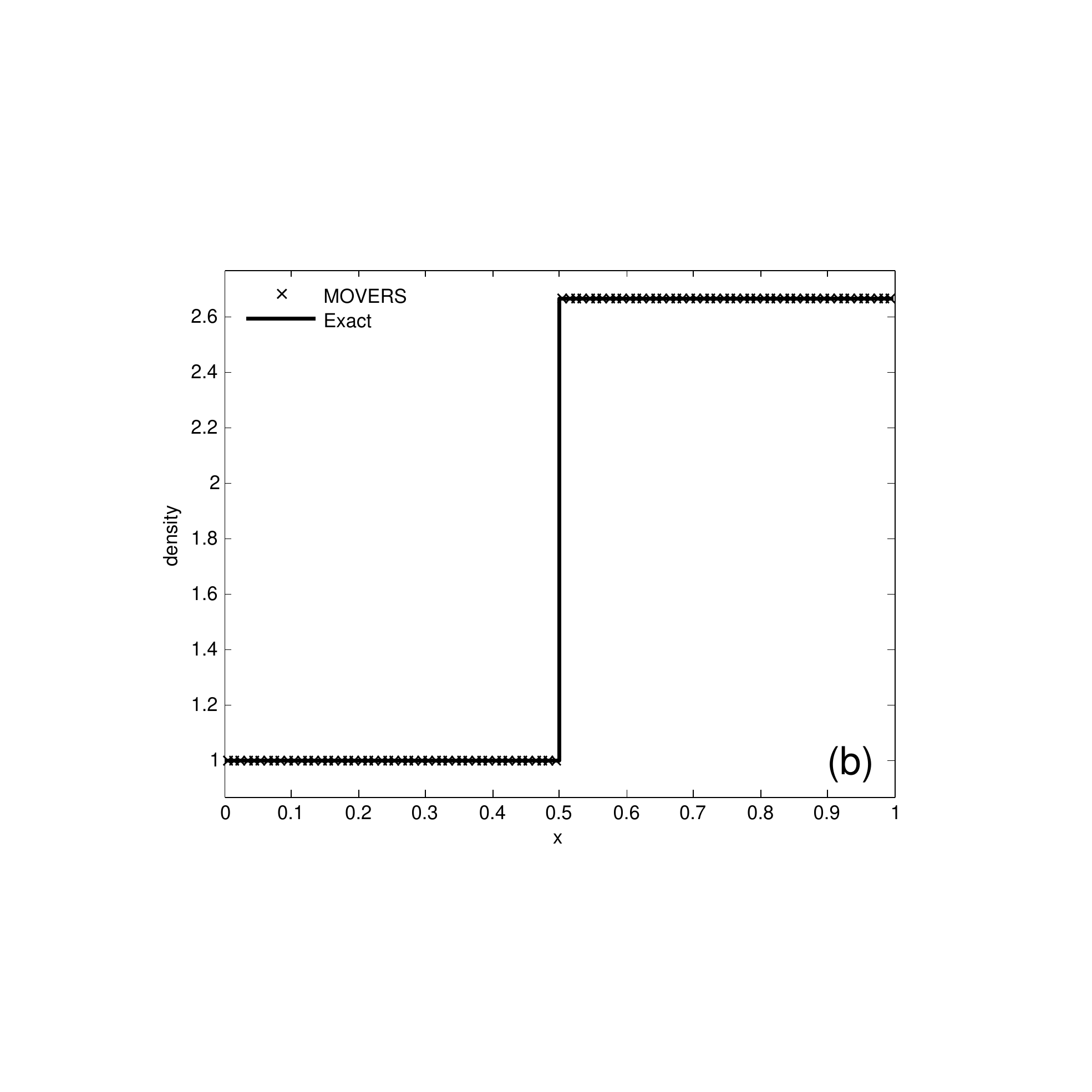}}
 \caption{Steady Contact Discontituity (a) and Steady Shock (b) results with MOVERS using 100 grid points}
  \label{stdM}
\end{figure}

\begin{figure}[h!] 
 \subfigure{\includegraphics[trim=75.0 125.0 100.0 131.0, clip, width=0.5\textwidth]{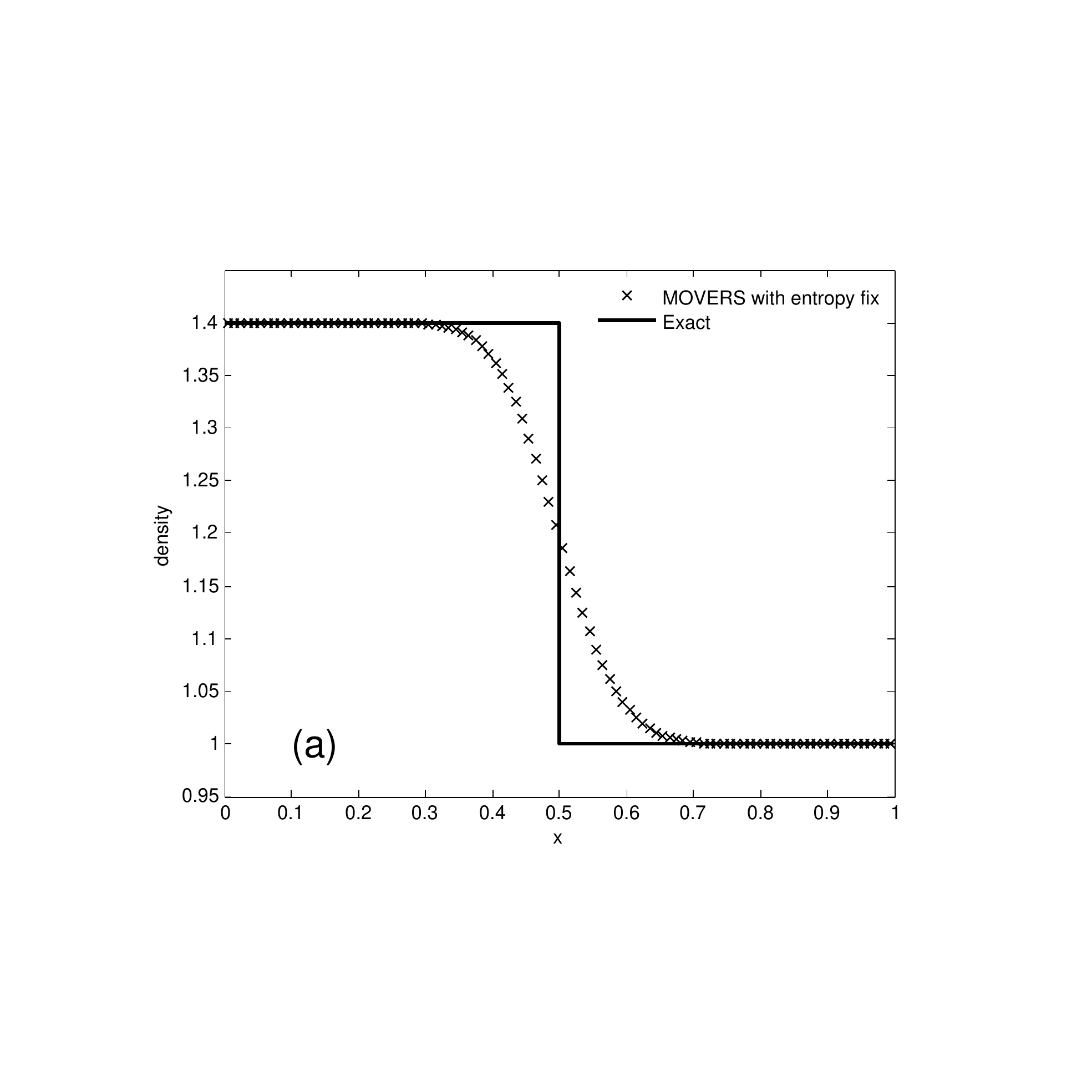}}
 \subfigure{\includegraphics[trim=75.0 125.0 100.0 131.0, clip, width=0.5\textwidth]{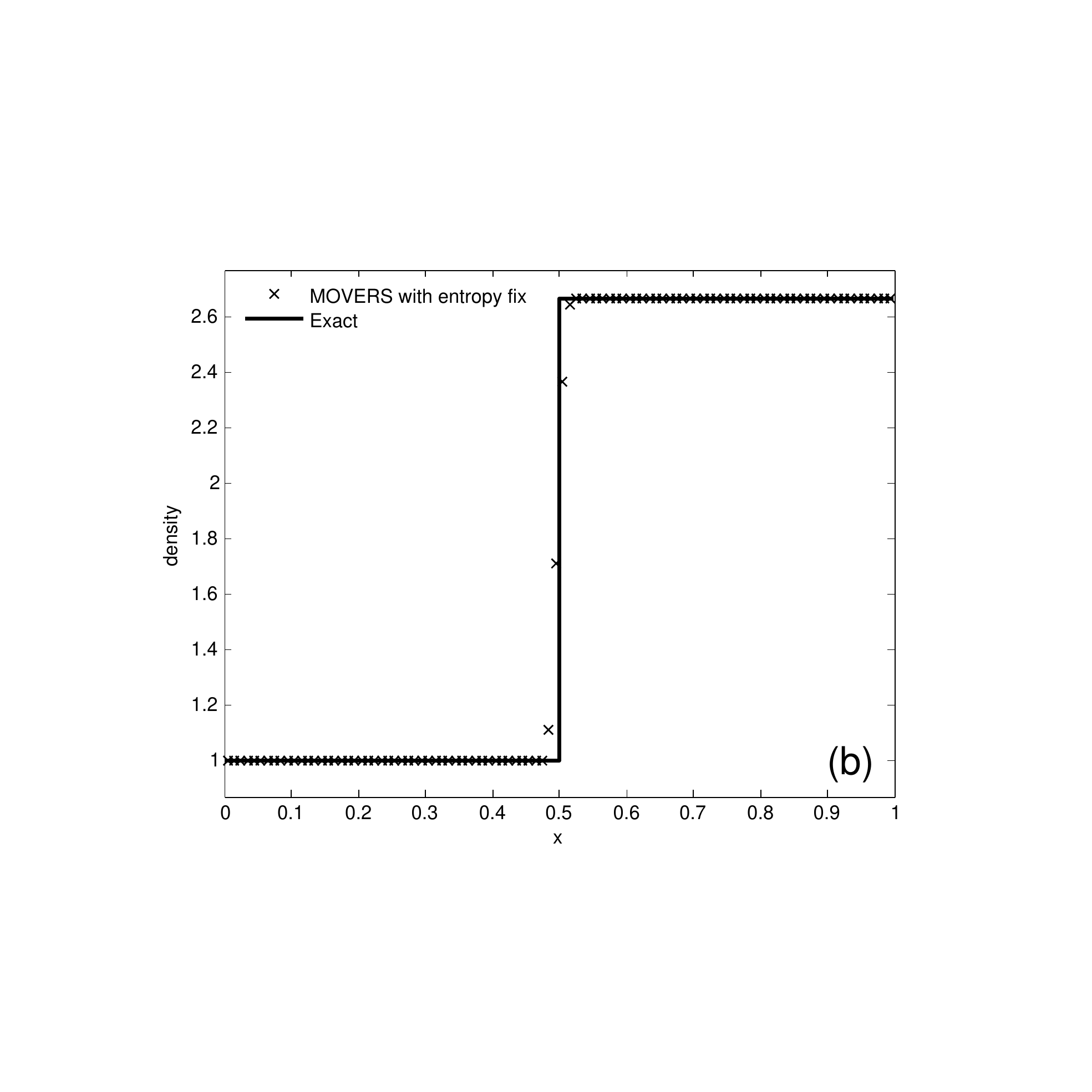}}
 \caption{Steady Contact Discontinuity (a) and Steady Shock (b) results with MOVERS and an entropy fix using 100 grid points}
 \label{stdMwef}
\end{figure}

\begin{figure}[h!]
  \subfigure{\includegraphics[trim=75.0 128.0 100.0 131.0, clip, width=0.5\textwidth]{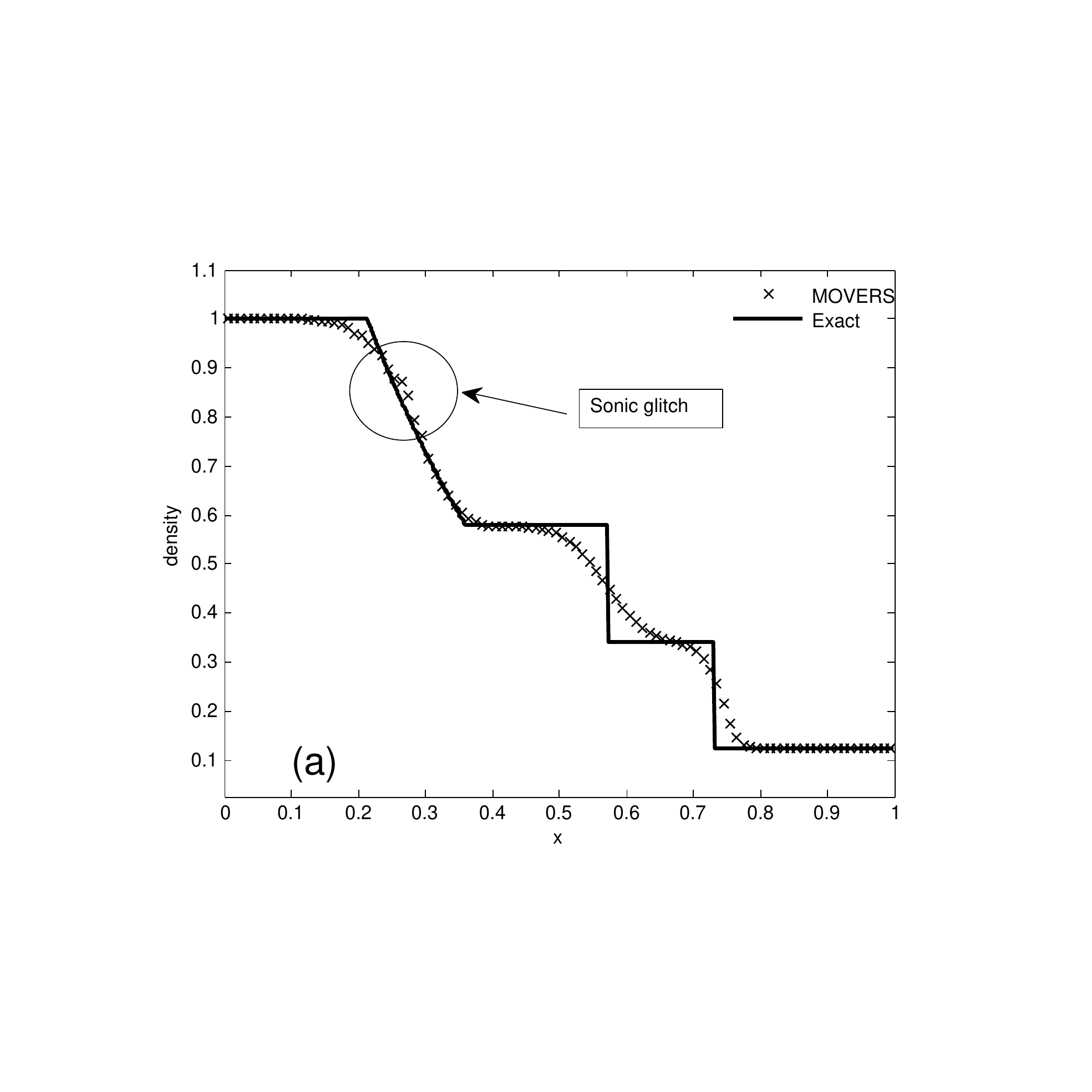}}
  \subfigure{\includegraphics[trim=75.0 128.0 100.0 131.0, clip, width=0.5\textwidth]{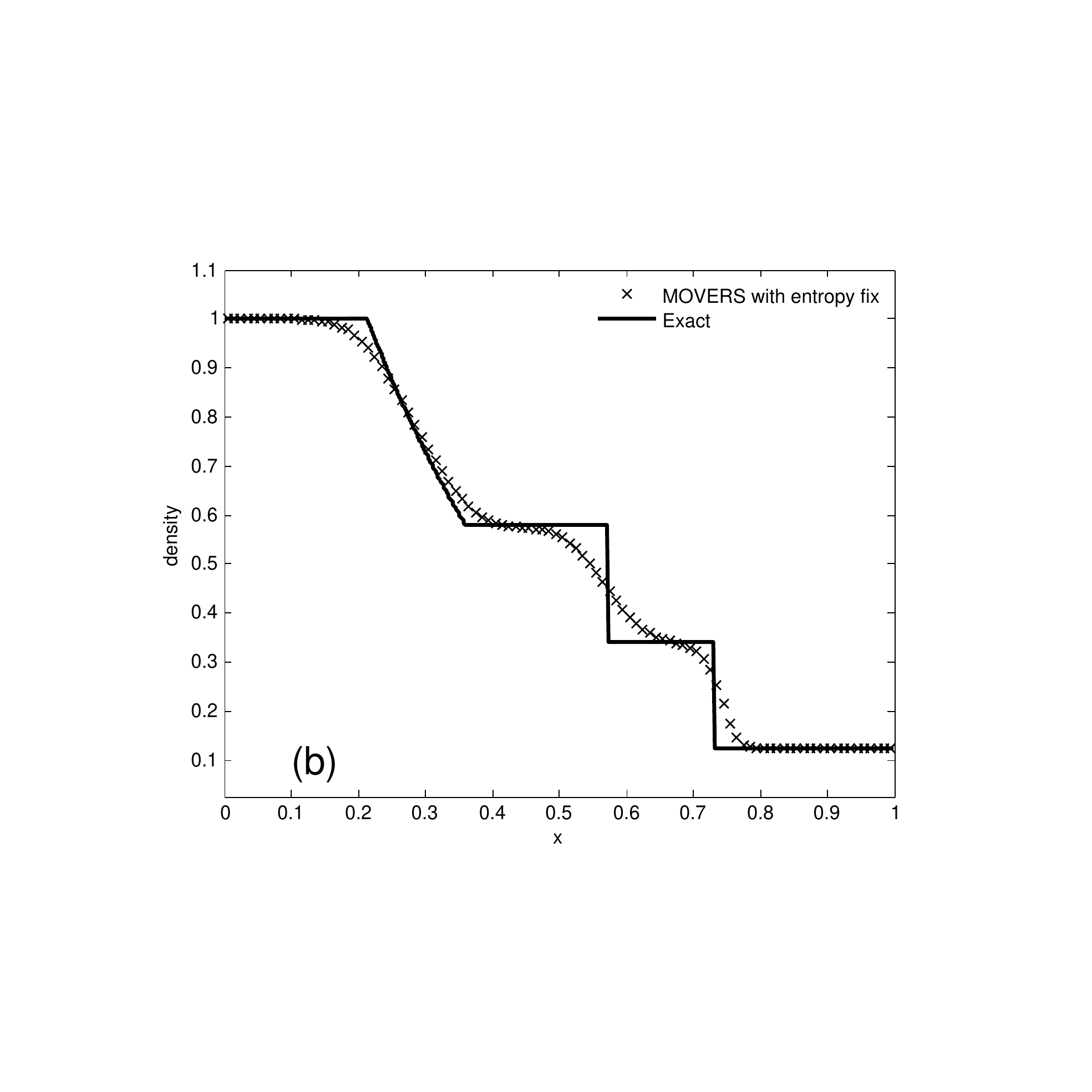}}
  \caption{First order result with MOVERS: (a) without and (b) with entropy fix for modified shock tube problem with sonic point using 100 grid points}
   \label{sonicglitch}
\end{figure}

The MOVERS is tested on steady contact and steady shock test cases to demonstrate its ability to capture them exactly. Details of these test cases are explained in the Section \ref{SecStdConShk}. As shown in the Fig (\ref{stdM}a), MOVERS captures steady the contact discontinuity exactly. The steady shock wave problem is taken from \cite{zhang2007new}, as shown in the  Fig (\ref{stdM}b) and MOVERS captures the steady shock also exactly without any numerical diffusion.  Modified shock tube problem is one of the test problems used to assess the numerical method for its entropy stability property. The details of the problem are given in the section \ref{sodshock}. MOVERS when tested on this problem, gives a sonic glitch in the expansion region as shown in the Fig (\ref{sonicglitch}a). MOVERS being an  accurate scheme, has less numerical diffusion, and thus requires an entropy fix to avoid non-smoothness in the expansion regions. The details of the entropy fix are explained in the paper \cite{jaisankar2009central}. It is worth noting here that Roe scheme \cite{roe1981approximate} gives an unphysical expansion shock, whereas MOVERS produces a milder sonic glitch. By using an entropy fix the problem of sonic glitch is eliminated, which is evident from Fig (\ref{sonicglitch}b).  With the use of an entropy fix, although sonic glitch problem is overcome, the ability of MOVERS to capture steady discontinuities exactly is lost.  As demonstrated in the Fig. (\ref{stdMwef}, the steady discontinuities are then diffused.  Motivated by the desire to overcome this problem, Jaisankar and Raghurama Rao \cite{jaisankar2009central} proposed a hybrid scheme named MOVERS-L, which uses the diffusion based on Rusanov scheme \cite{rusanov1962calculation} in smooth regions of the flow and MOVERS at the discontinuities, based on a limiter function which is discussed in the next section.

\subsection{MOVERS-L}
\label{SecMoversL}
The motivation for this version of the MOVERS is to retain the property of exact discontinuity capturing in steady state while avoiding the entropy fix.  The Rusanov (LLF) scheme \cite{rusanov1962calculation} is a simple central solver  which satisfies entropy condition and avoids expansion shocks.  The coefficient of numerical diffusion in the LLF method is fixed as the maxium of the eigenvalues of the flux Jacobian matrix, from the left and rights states of the cell interface.  In MOVERS-L, the numerical diffusion based on LLF scheme is used in the smooth regions of the flow and the numerical diffusion based on R-H condition as in MOVERS is used at the large gradients.  The switching over between smooth regions and large gradients of the flow is done using a minmod limiter \cite{yee1987construction}. The coefficient of numerical diffusion in this method, at the cell interface, is given below.

\begin{equation} \label{moversle}
 \alpha_{j+\frac{1}{2}}=\alpha_{MOVERS}-\phi(r) (\alpha_{MOVERS}-\alpha_{LLF})
\end{equation}
where $\alpha_{MOVERS}$ is the diffusion given by the expressions \eqnref{rhc} to \eqnref{mineig}, and $\alpha_{LLF}$ is the maximum spectral radius of the flux Jacobian matrix from the left and right states of the cell interface. $\phi(r)$ is the minmod flux limiter which acts as a switch and distinguishes between smooth regions and high gradients of the flow.

\begin{equation} \label{minmod}
 \phi_{j+\frac{1}{2}} \left( r_{j+\frac{1}{2}} \right) = minmod(1, r_{j}^+, r_{j+1}^- )
\end{equation}

\begin{equation} \label{rj+}
  r_{j}^+ =\frac{U_j-U_{j-1}}{U_{j+1}-U_j} ~~\mbox{and}~~ r_{j}^- =\frac{1}{r_{j}^+}
\end{equation}
In the above expression of $r_{j}^+$, as the denominator is the difference in conserved variables, when $U_{j+1}\approx U_j$ the value of the expression can lead to large values of $r_{j}^{+}$.  To avoid the numerical overflow, the following definition of $r_{j}^+$ is used.

\begin{equation} \label{rj+mod}
  r_{j}^+ = sign(U_{j+1}-U_j) \frac{U_j-U_{j-1}}{\delta} ~~\mbox{if}~~ |U_{j+1}-U_j|<\delta
\end{equation}
where $\delta$ is small number.
The multivariable $minmod$ function is defined as follows \cite{Kurganov_Tadmor}.  

\begin{equation}
 minmod(x_{1}, x_{2}, \ldots)=\begin{cases}
    min_{j}\left\{x_{j}\right\} \ \textrm{if} \ x_{j} > 0 \ \forall j \\[2mm] 
    max_{j}\left\{x_{j}\right\} \ \textrm{if} \ x_{j} < 0 \ \forall j \\[2mm] 
    0 \ \textrm{otherwise} 
    \end{cases}
\end{equation} 

\begin{figure}[h!] 
 \subfigure{\includegraphics[trim=75.0 125.0 100.0 131.0, clip, width=0.5\textwidth]{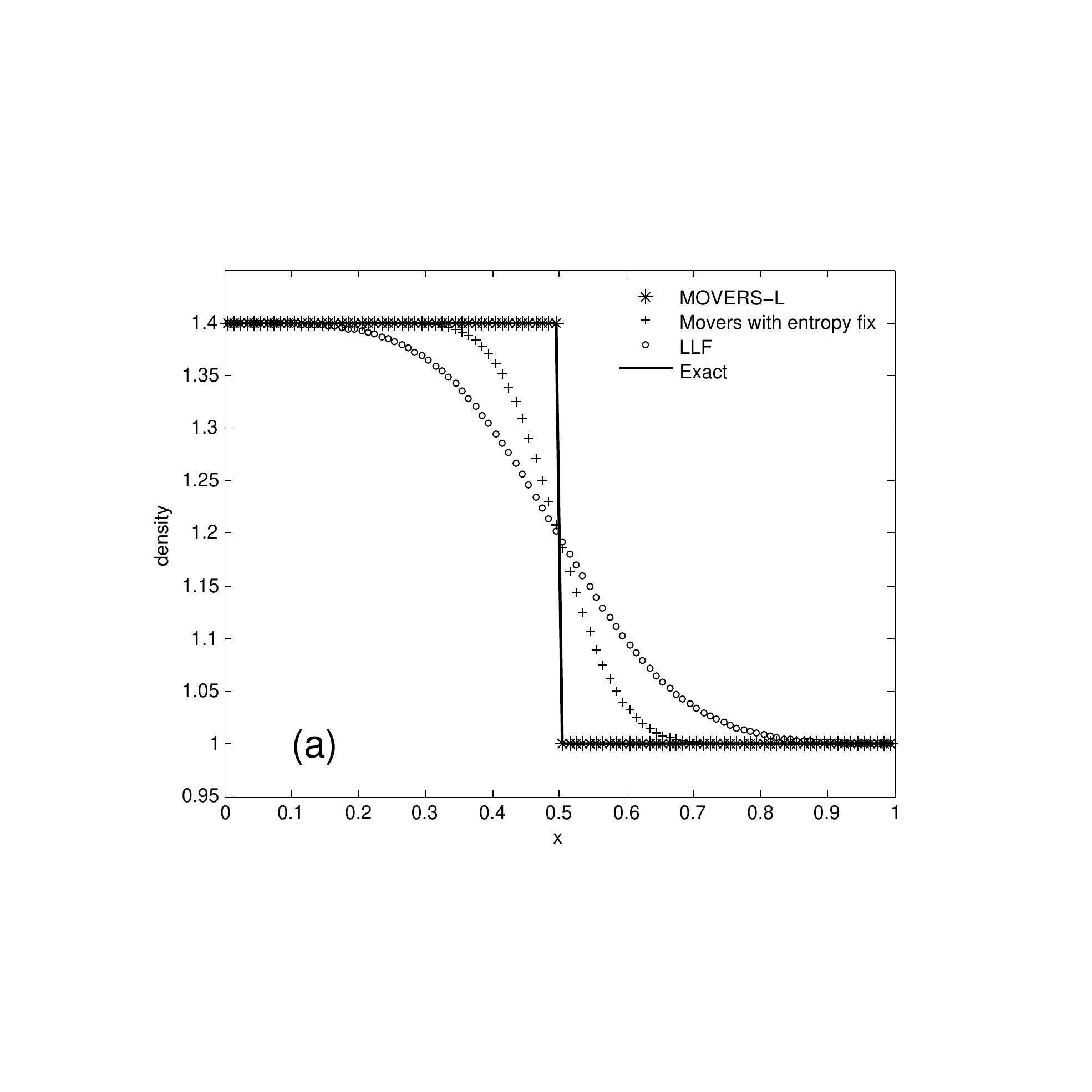}}
 \subfigure{\includegraphics[trim=75.0 125.0 100.0 131.0, clip, width=0.5\textwidth]{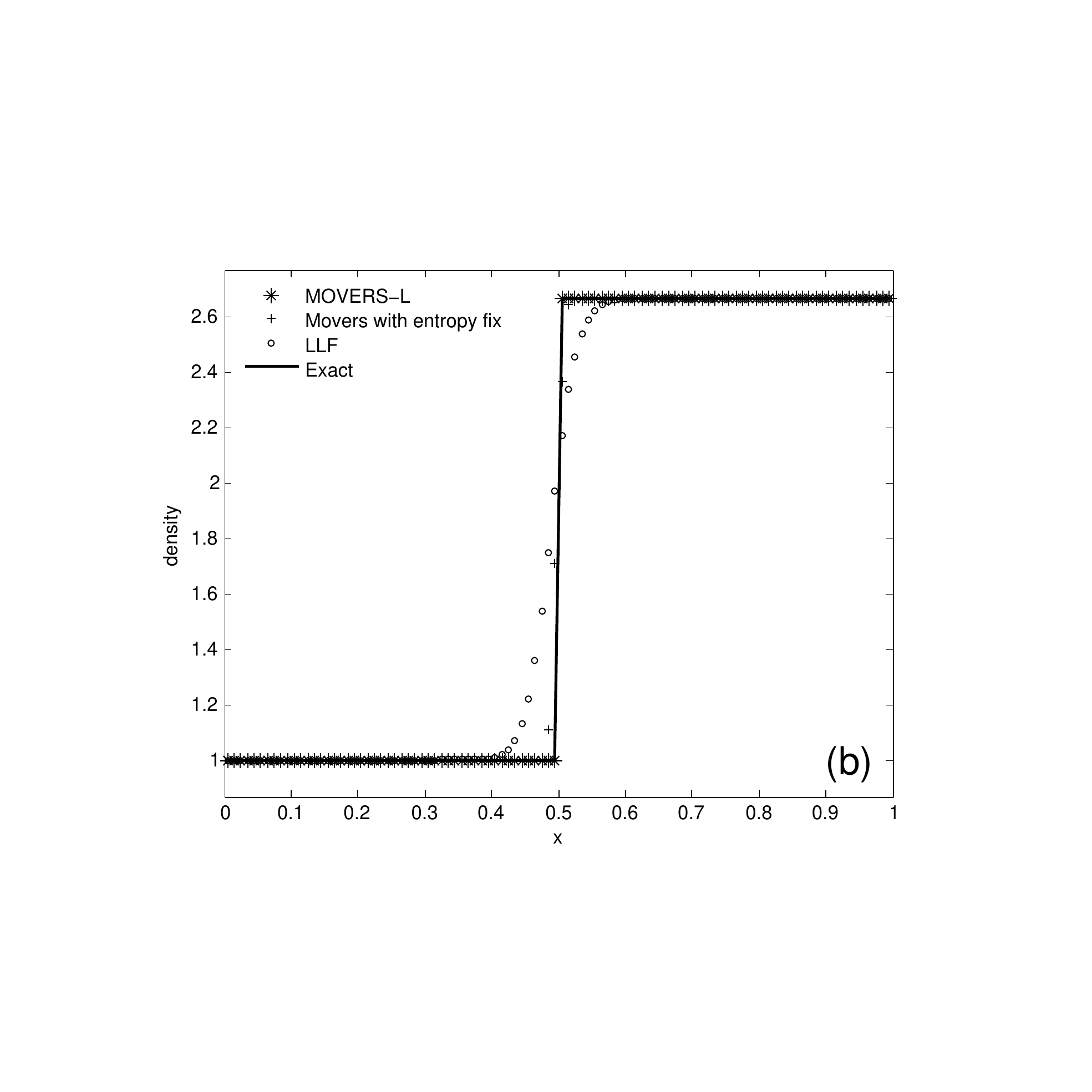}}
 \caption{Comparison of MOVERS-L results with MOVERS and LLF for Steady Contact Discontinuity (a) and Steady Shock (b) }
 \label{stdMoversL}
\end{figure}

\begin{figure}[h!] 
\begin{center}
 \includegraphics[trim=85.0 130.0 100.0 131.0, clip, width=0.8\textwidth]{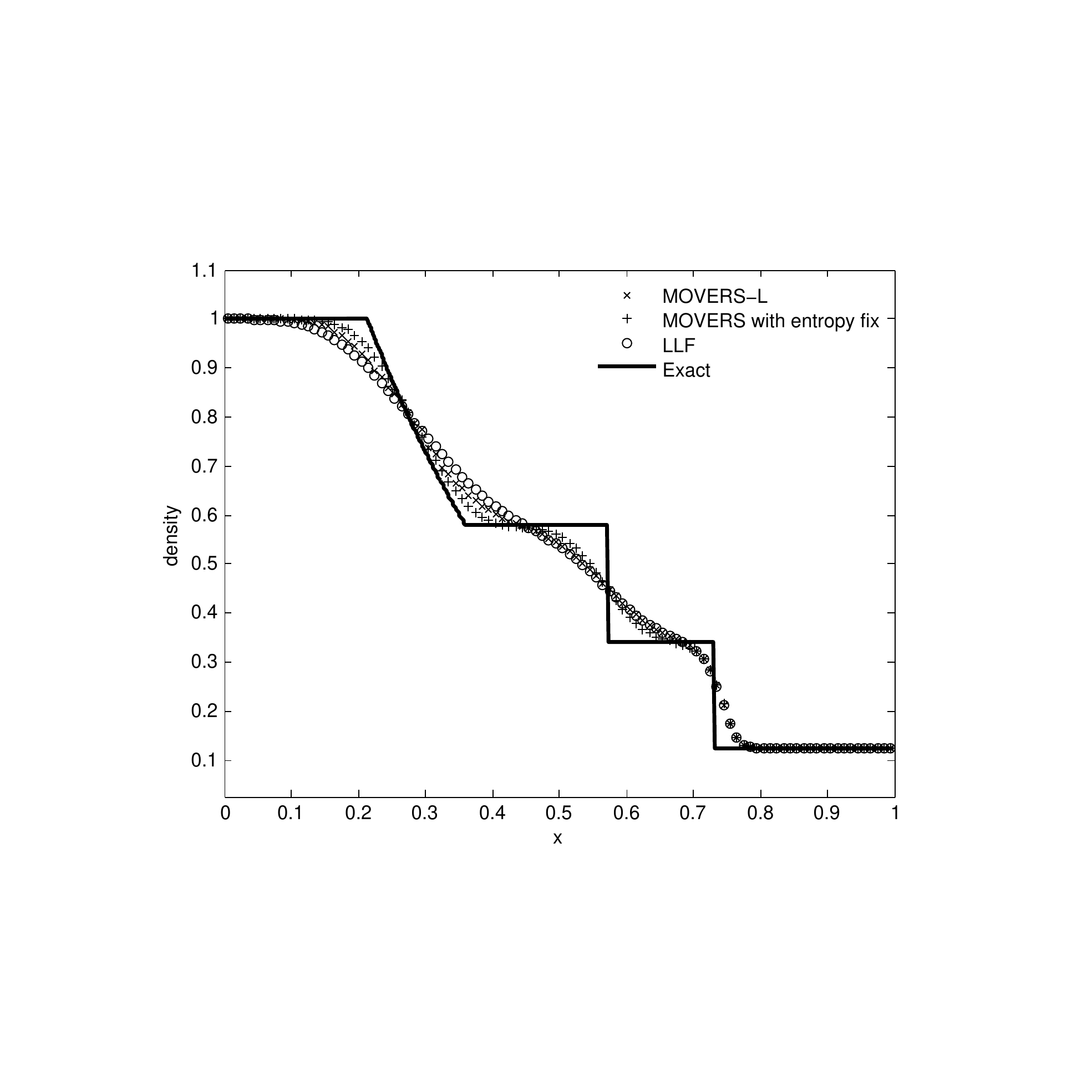}
 \caption{Comparison of MOVERS-LE result with LLF and MOVERS (with entropy fix) for modified shock tube problem with sonic point using 100 grid points}
 \label{FigMoversLSonic}
\end{center}
\end{figure}

As expected, MOVERS-L captures steady discontinuities exactly which is shown in the Fig.~(\ref{stdMoversL}). The results are also compared with LLF and MOVERS with entropy fix.  This scheme is tested on modified shock tube problem with a sonic point in the expansion fan to test the entropy stability property, as shown in Fig (\ref{FigMoversLSonic}).  As expected, MOVERS-L avoids generating an expansion shock.  Although MOVERS-L captures steady discontinuities exactly and yet avoids using an entropy fix, it is diffusive in smooth regions as is evident from Fig (\ref{FigMoversLSonic})  since we are using LLF based diffusion in the smooth regions of the flow.   Developing a numerical method which captures steady discontinuities exactly, is accurate in smooth regions and yet is entropy stable is the motivation for the present work.  

\subsection{Entropy conservation equation based numerical diffusion}\label{SecEntrConser}
The main focus of this work is to enforce entropy stability in smooth regions, retaining sharp resolution of discontinuities.  The expansion regions with sonic points require finite numerical diffusion to avoid any expansion shocks, while shocks need less diffusion to resolve them accurately.  Ensuring both these conflicting requirements in a numerical method is a non-trivial task.  The basic idea here is to use the entropy conservation equation as a guide-line to fix the coefficient of numerical viscosity.  Considering the second law of thermodynamics together with the conservation of mass, momentum and energy equations, we can rewrite the 1-D Euler equations using an entropy formulation as follows 
\cite{leveque1992numerical, Whitham, Naterer_Camberos}.  
\begin{eqnarray} 
\frac{\partial \rho}{\partial t} + \frac{\partial \left( \rho u \right)}{\partial x} = 0 \\ 
\frac{\partial \left( \rho u \right)}{\partial t} + 
\frac{\partial \left( p + \rho u^{2} \right)}{\partial x} = 0 \\ 
\frac{\partial \left(\rho S\right)}{\partial t} 
+ \frac{\partial \left(\rho u S\right)}{\partial x} = 0   
\end{eqnarray} 
where entropy $S$ is defined by 
\begin{equation} 
S = C_{V} \log \left( \frac{p}{\rho^{\gamma}}\right) + constant
\end{equation} 
Note that the energy conservation equation is replaced by an entropy conservation equation after a little bit of algebraic manipulation of the equations.   This formulation is valid only for smooth regions, without discontinuities.  For flows with discontinuities, the entropy conservation equation has to be replaced by an inequality as follows.  
\begin{equation} \label{entropyineq}
\frac{\partial{(\rho S)}}{\partial{t}} +  \frac{\partial{(\rho u S)}}{\partial{x}} \geq 0 \linebreak
\end{equation}  
To keep the entropy bounded, mathematicians define the entropy with a negative sign so that we obtain the following equation.  
\begin{equation} \label{entropyineq2}
\frac{\partial{(\rho S)}}{\partial{t}} +  \frac{\partial{(\rho u S)}}{\partial{x}} \leq 0 \linebreak
\end{equation}  
In the absence of discontinuities (for smooth regions of the flow), the equality holds and the entropy conservation equation is given as
\begin{equation} \label{entropy_conservation} 
\frac{\partial{(\rho S)}}{\partial{t}} +  \frac{\partial{(\rho u S)}}{\partial{x}} = 0  \linebreak
\end{equation}
Tadmor\cite{tadmor2003entropy} has analyzed entropy stability and entropy conservation properties of numerical methods.  He argues that, for any scheme to be entropy stable, it should satisfy additional inequality \eqnref{entropyineq2}.  The basic strategy of the present research work is to use the above entropy conservation equation to get some information about how to fix the coefficient of numerical diffusion in finite volume method in smooth regions of the flow.  Our approach is much simpler than that of Tadmor, while Tadmor's approach is mathematically rigorous.  

\begin{figure}[h!] 
 \subfigure{\includegraphics[trim=75.0 125.0 100.0 131.0, clip, width=0.5\textwidth]{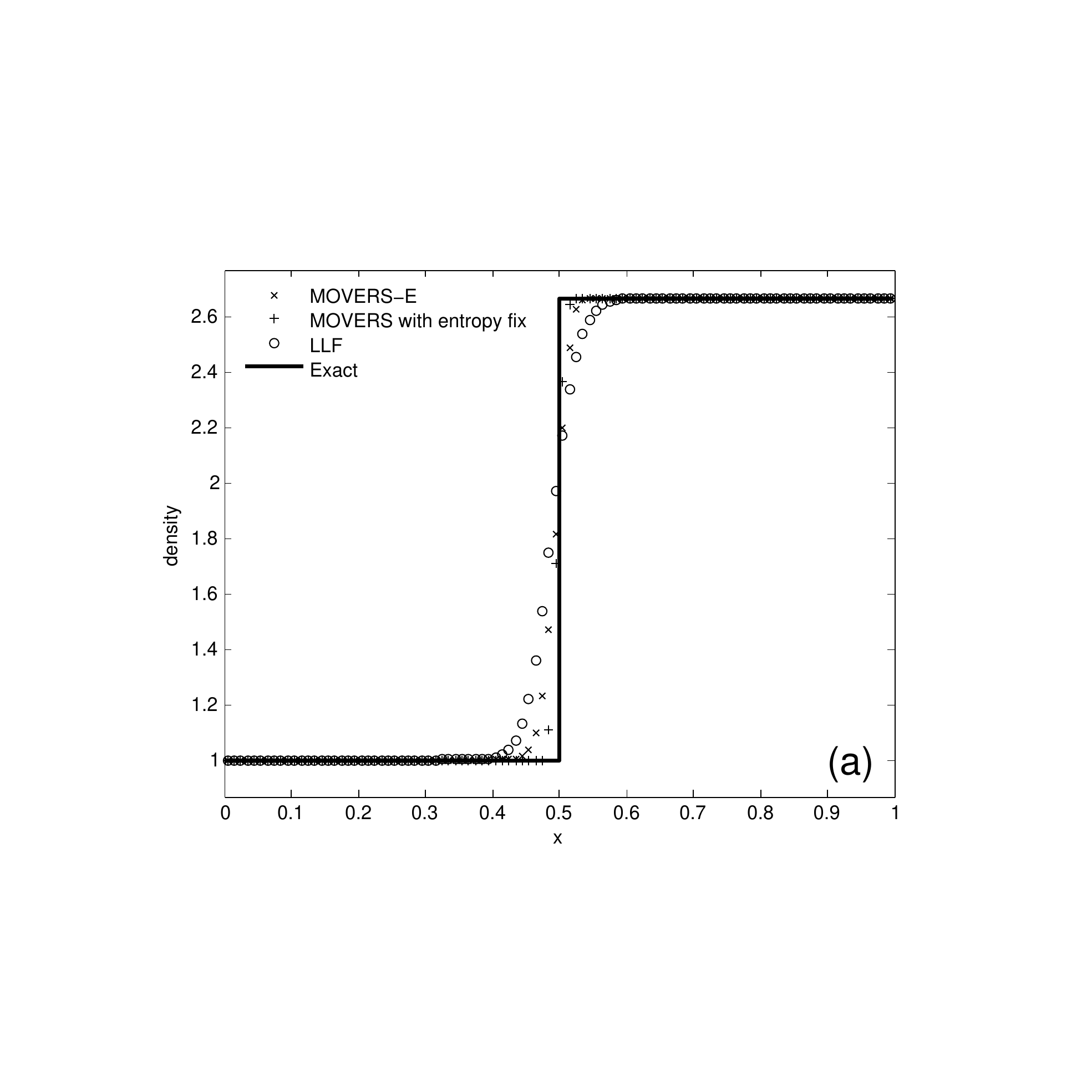}}
 \subfigure{\includegraphics[trim=75.0 125.0 100.0 131.0, clip, width=0.5\textwidth]{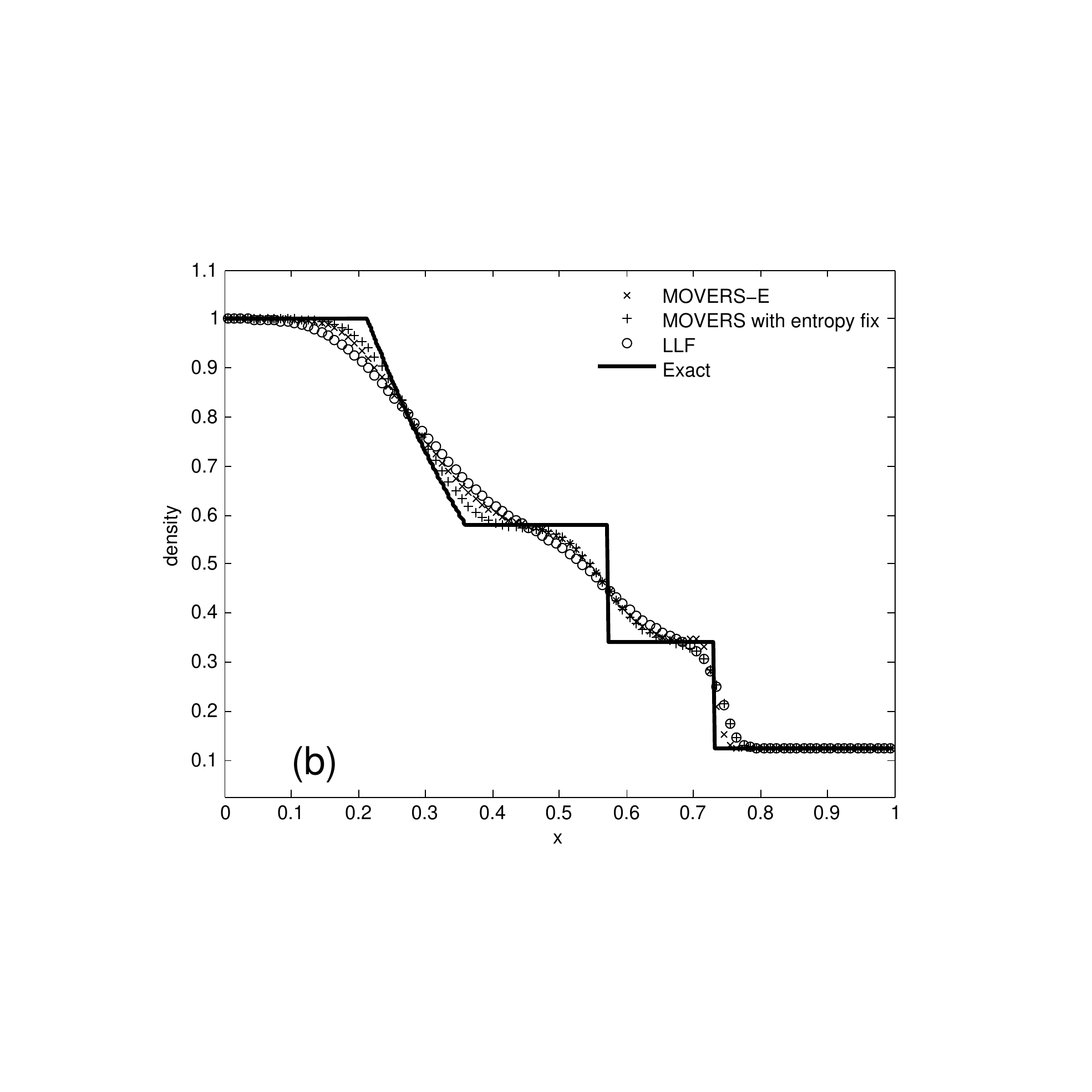}}
 \caption{Comparison of MOVERS-E results with MOVERS and LLF: (a) Steady Contact Discontinuity, (b) Modified shock  tube with sonic point }
 \label{stdMoversE}
\end{figure}

A very efficient and yet simple strategy is to apply the LLF method for the entropy conservation equation to obtain the numerical diffusion and to simply use the same numerical diffusion for discretizing the Euler equations, using a finite volume method.  Since the entropy conservation equation is applicable for smooth regions and moreover it is used in place of the energy conservation equation, the hope is that discretization will yield the required amount of numerical diffusion to yield a stable scheme for all the conservation equations in the smooth regions of the flow.   It is also expected to yield an entropy stable scheme for Euler equations.  Eventually, the plan is to couple this strategy to the strategy of MOVERS being used near the large gradients, thereby retaining the exact disctoninuity capturing in a simple central discretization framework.  Since the numerical diffusion in the LLF method is simply the local maximum of the eigenvalues of the flux Jacobian matrix for any system of equations, we first determine the wave-speed for the entropy conservation equation, which is given by 
\begin{equation}
a_{ES}= \frac{\partial{(\rho u S)}}{\partial{(\rho S)}} \approx u  \linebreak
\end{equation}
Thus, the numerical diffusion rquired for applying the LLF method to the entropy conservation equation will be
\begin{equation} \label{EqnEntrStab}
 a_{LLF- ES} = max[|u_L|,|u_R|] \linebreak
\end{equation}
In $a_{ES}$, the subscript $ES$ refers to the phrase {\em Entropy Stable (ES)},  where L and R denote the left and right states of a cell interface.  Since this numerical diffusion is from a stable discretization of the entropy conservation equation for smooth regions, we can use the same numerical diffusion for discretizing Euler equations, with the expectation of obtaining a stable numerical method for smooth regions. This method is termed as MOVERS-E.
This scheme is tested on steady shock problem and shock tube problem with sonic point. As shown in the Fig (\ref{stdMoversE}b) the scheme has no sonic point problem, and yet is accurate. But the scheme diffuses the steady shocks Fig (\ref{stdMoversE}a). By following the strategy of MOVERS-L, we device a hybrid scheme named MOVERS-LE (L refers to limiter and E refers to entropy stability) as described in the next section.

\subsection{Entropy stable hybrid solver} \label{SecEntrStab}

The special discretization described in  section \ref{SecEntrConser} is an accurate and stable numerical scheme for only smooth flows as it diffuses the discontinuities. Near discontinuities, we can use the numerical diffusion based on Rankine-Hugoniot condition, as in MOVERS for sharp resolution of the high gradients. This hybrid scheme, named MOVERS-LE, can be obtained by utilizing the following coefficient of numerical diffusion at the interface.

\begin{equation}
 \alpha_{MOVERS-LE, j+\frac{1}{2}}=\alpha_{MOVERS, j+\frac{1}{2}} 
+ \phi (r_{j+\frac{1}{2}})\left[\alpha_{LLF-ES, j+\frac{1}{2}} 
- \alpha_{MOVERS, j+\frac{1}{2}}\right]  
\end{equation} 
where, $\alpha_{LLF-ES, j+\frac{1}{2}}$ is given by \eqnref{EqnEntrStab} and 
$\alpha_{MOVERS}$ is the R-H condition based numerical diffusion, as explained in Section \ref{SecMoversL}. The results of this hybrid scheme, MOVERS-LE, is presented in the next section for various test cases.  In two dimensional case, a finite volume method is used with the above strategy being implemented across any cell interface in the direction normal to the interface.  

\section{Results and discussions}
 
\subsection{1D-Euler equation results}
MOVERS-LE is tested on following 1D-Euler test cases, which are taken from Toro~\cite{toro1999riemann} to assess the numerical method for accuracy and robustness.  

\subsubsection{Steady shock and contact discontinuity} \label{SecStdConShk}
\begin{figure}[h!] 
 \subfigure{\includegraphics[trim=75.0 125.0 100.0 131.0, clip, width=0.5\textwidth]{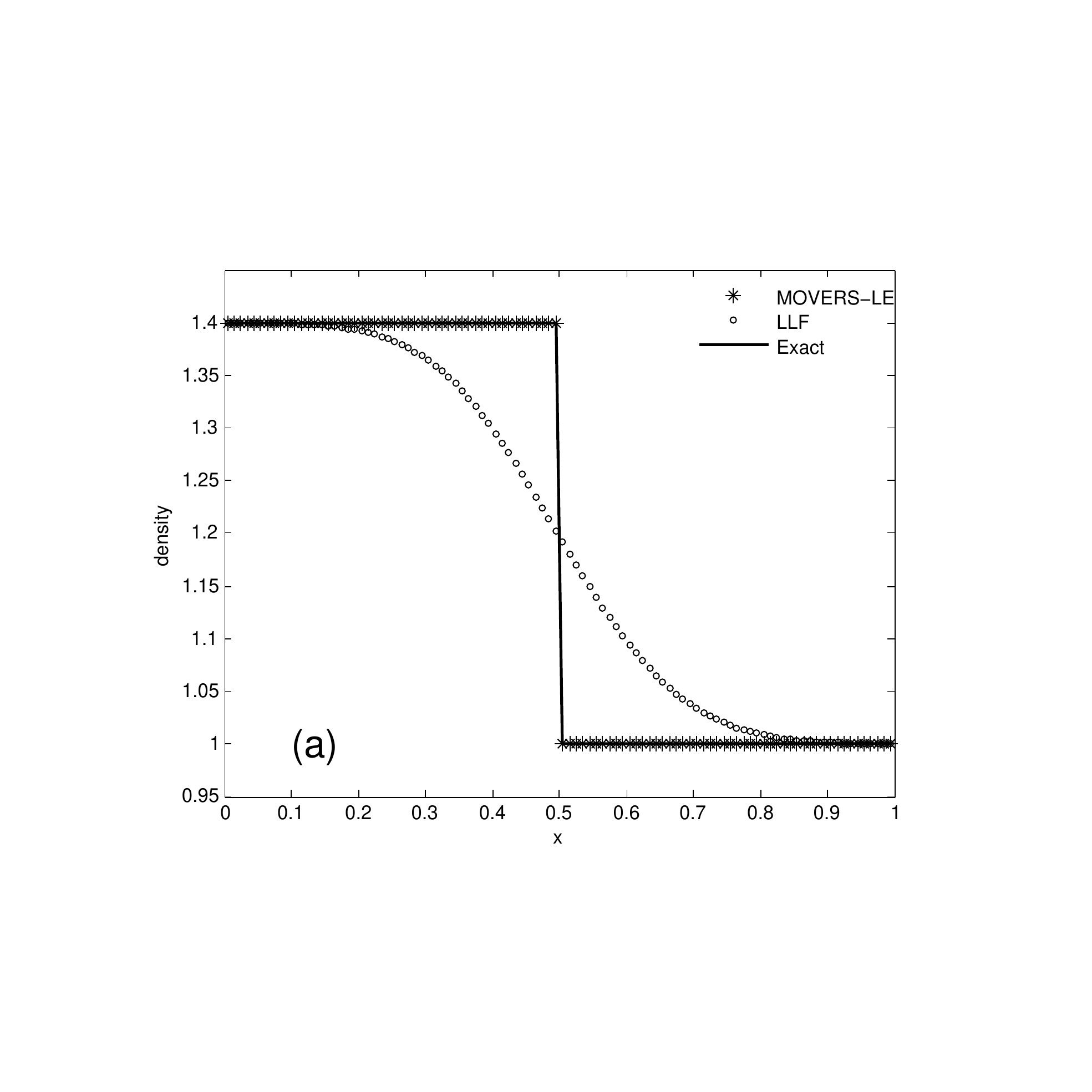}}
 \subfigure{\includegraphics[trim=75.0 125.0 100.0 131.0, clip, width=0.5\textwidth]{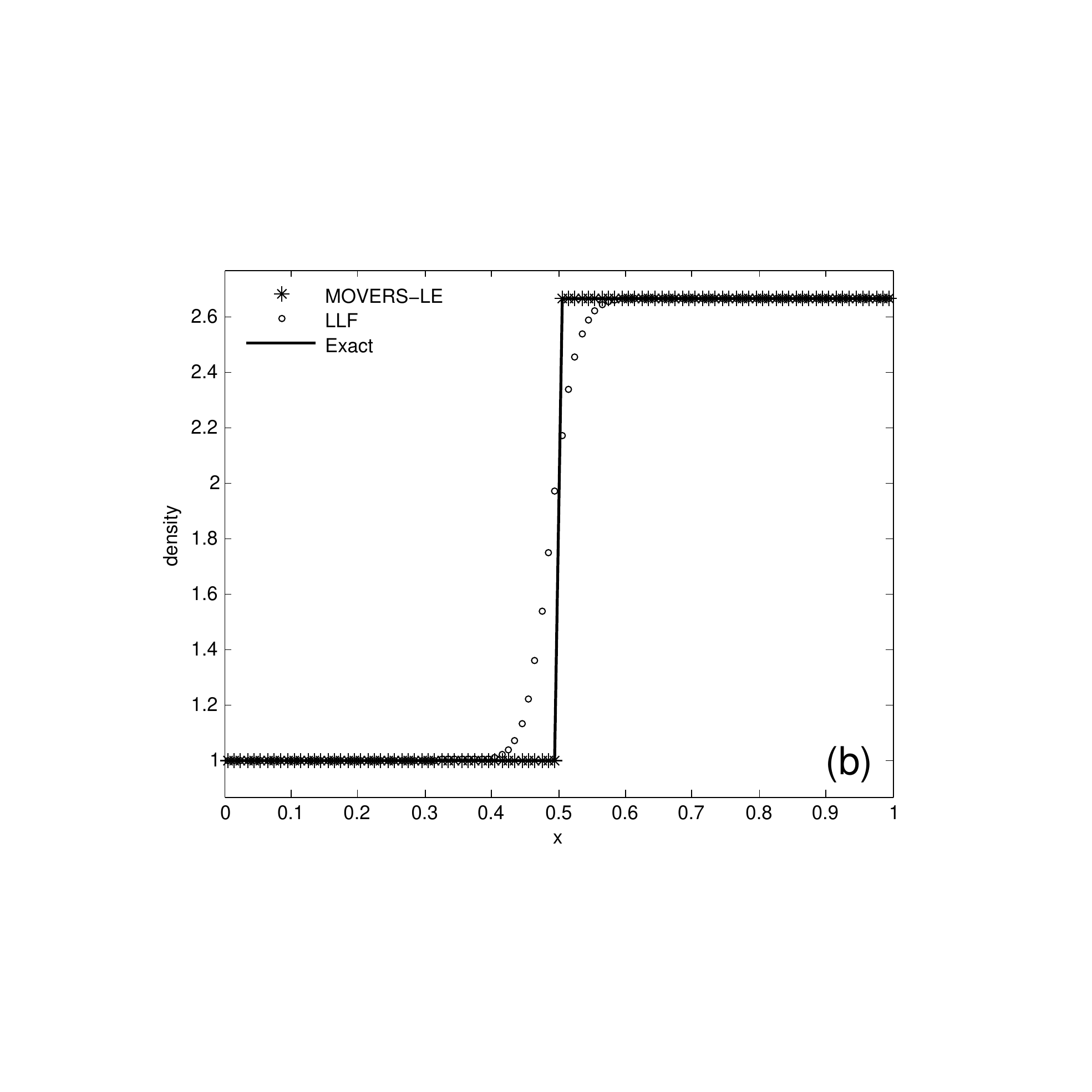}}
 \caption{Comparison of MOVERS-LE results with LLF for: (a) Steady Contact Discontinuity, (b) Steady Shock}
 \label{FigstdMoversLE}
\end{figure}

For the steady contact discontinuity, the initial condition are given as
\begin{equation*}
(\rho, u, p)_R=(1.4, 0.0, 1.0) ~~\mbox{and}~~ (\rho, u, p)_L=(1.0, 0.0, 1.0) 
\end{equation*}
The initial discontinuity is at $x_o=0.5$ and the final computation time is at $t=2.0$.  As shown the  Fig(\ref{FigstdMoversLE}a), MOVERS-LE captures this exactly. The steady shock problem is taken from Zhang and Shu \cite{zhang2007new}. This test case is built using R-H conditions.  MOVERS-LE captures steady shock also exactly without any numerical diffusion, as shown in the Fig (\ref{FigstdMoversLE}b).

\subsubsection{Sod's modified Shock tube} \label{sodshock}
\begin{figure}[h!] 
 \subfigure{\includegraphics[trim=75.0 125.0 100.0 131.0, clip, width=0.5\textwidth]{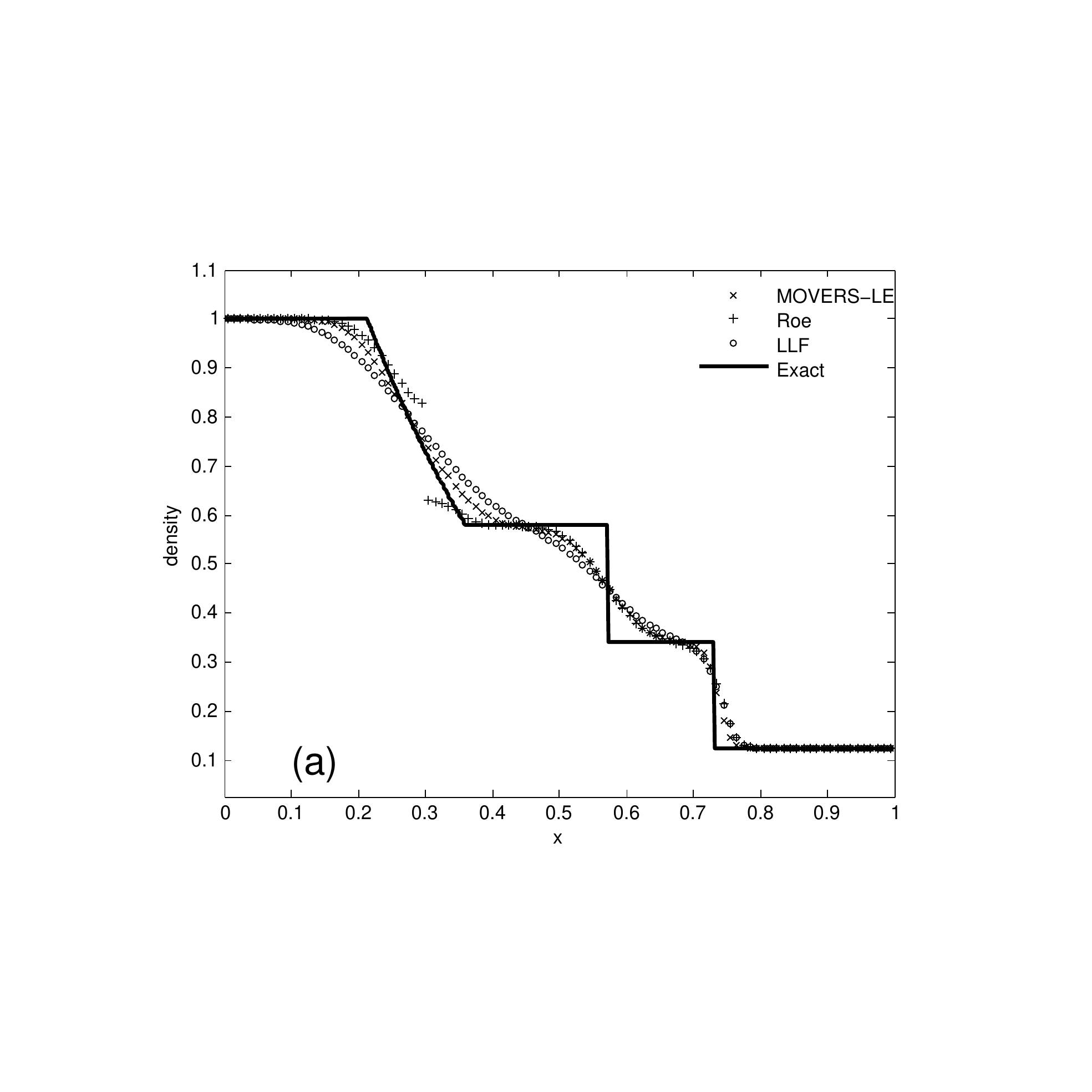}}
 \subfigure{\includegraphics[trim=75.0 125.0 100.0 131.0, clip, width=0.5\textwidth]{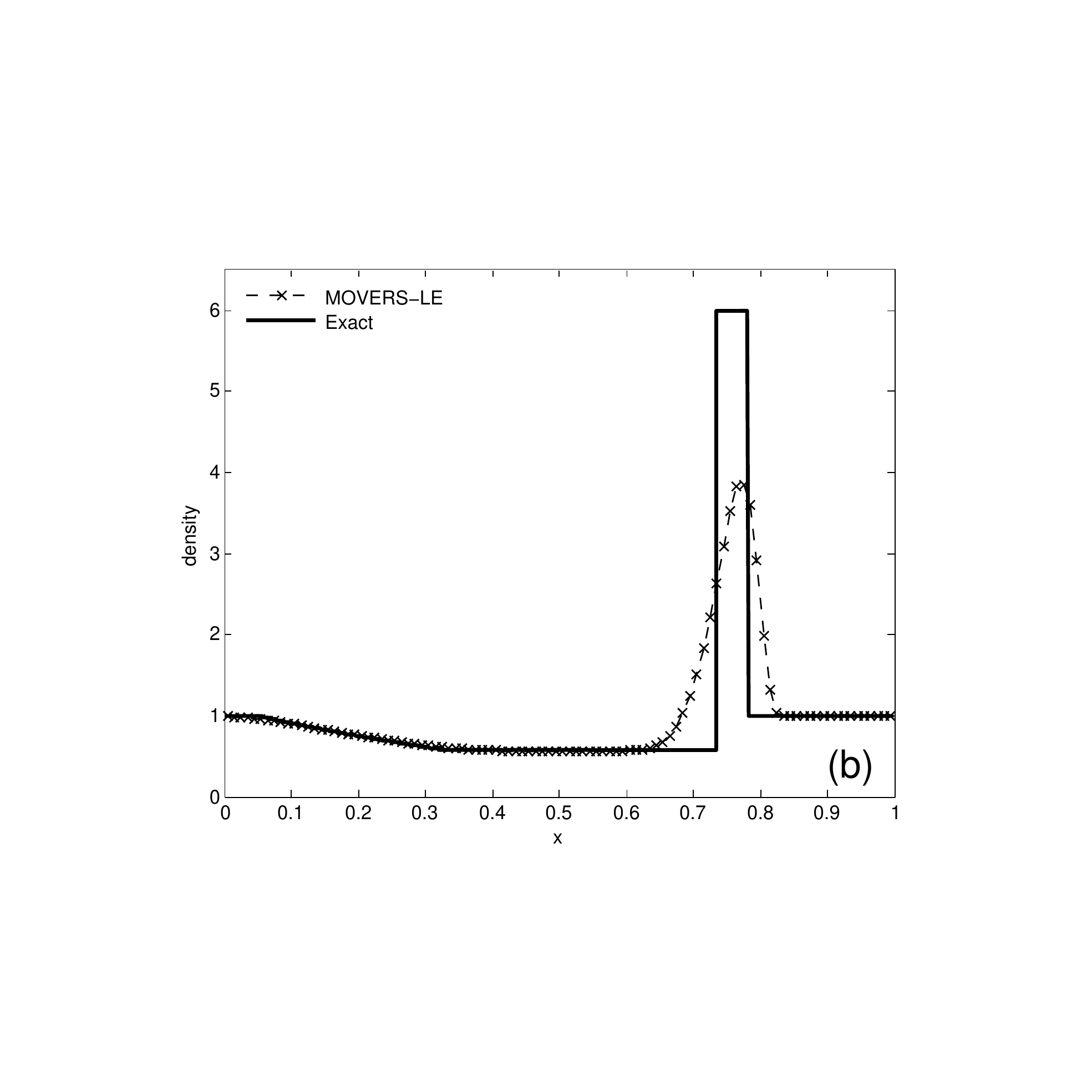}}
 \caption{(a) Comparison of MOVERS-LE result with LLF and Roe schemes for modified shock tube problem with sonic point using 100 grid points, (b) Strong Shock }
 \label{FigMoversLEProb13}
\end{figure}

This test case is the modified version of Sod's \cite{sod1978survey} test case. The solution of this problem consists of shock and contact discontinuities which are moving to the right and an expansion with sonic point to the left. This test problem is used to assess the numerical method for its entropy stability.  As shown in the Fig(\ref{FigMoversLEProb13}a), Roe scheme \cite{roe1981approximate} gives an unphysical expansion shock at the sonic point, whereas MOVERS-LE gives a smooth expansion profile.   

\subsection{Strong shock}
This test case is used to check the robustness of the numerical method as it involves a strong shock wave of shock Mach number $198$. The solution of this problem has a contact discontinuity, a strong shock and an expansion wave. The Fig. (\ref{FigMoversLEProb13}b) shows the solution of this problem using the proposed scheme, which demonstrates its robustness.

\subsection{Strong discontinuities}
\begin{figure}[h!] 
 \subfigure{\includegraphics[trim=85.0 125.0 100.0 131.0, clip, width=0.5\textwidth]{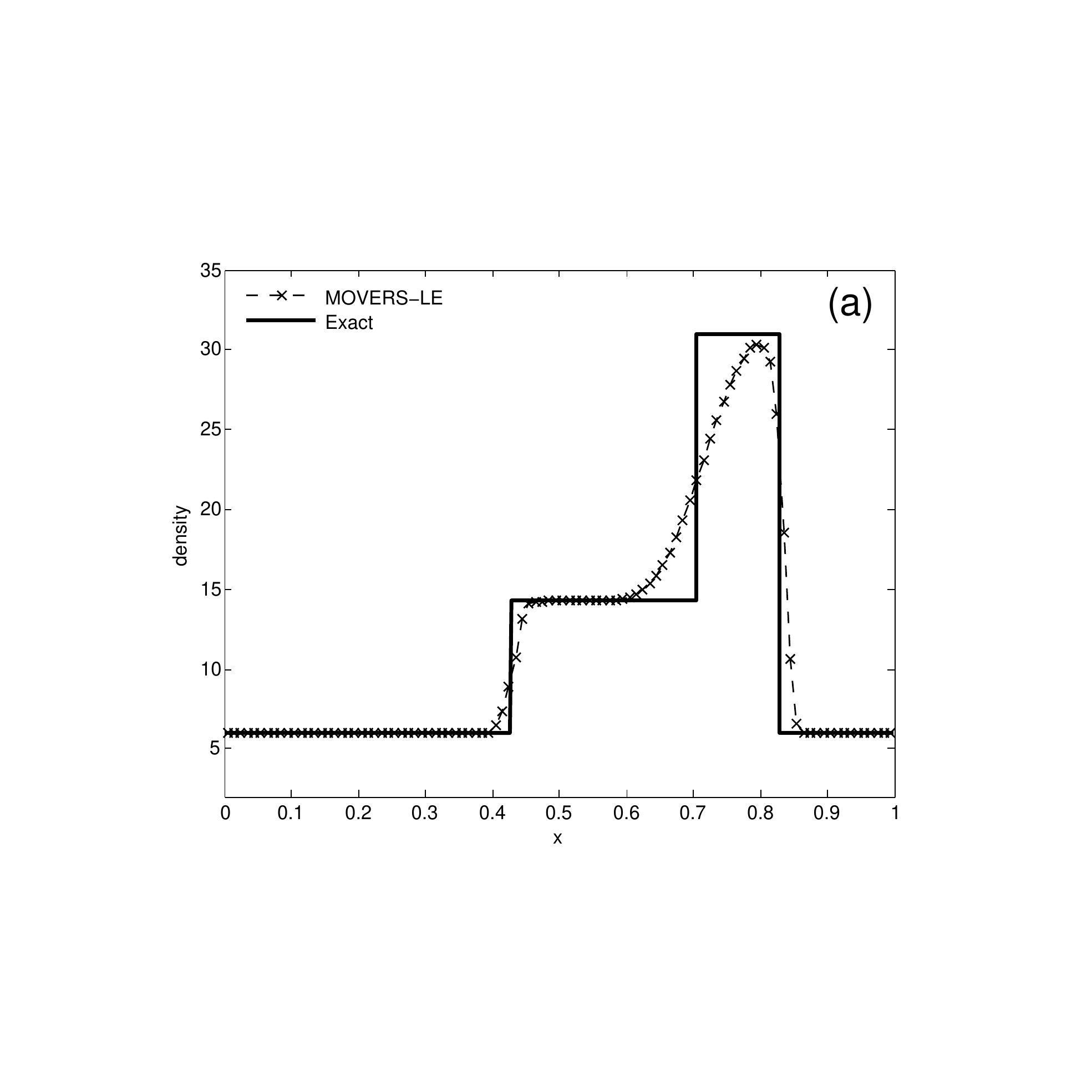}}
 \subfigure{\includegraphics[trim=85.0 125.0 100.0 131.0, clip, width=0.5\textwidth]{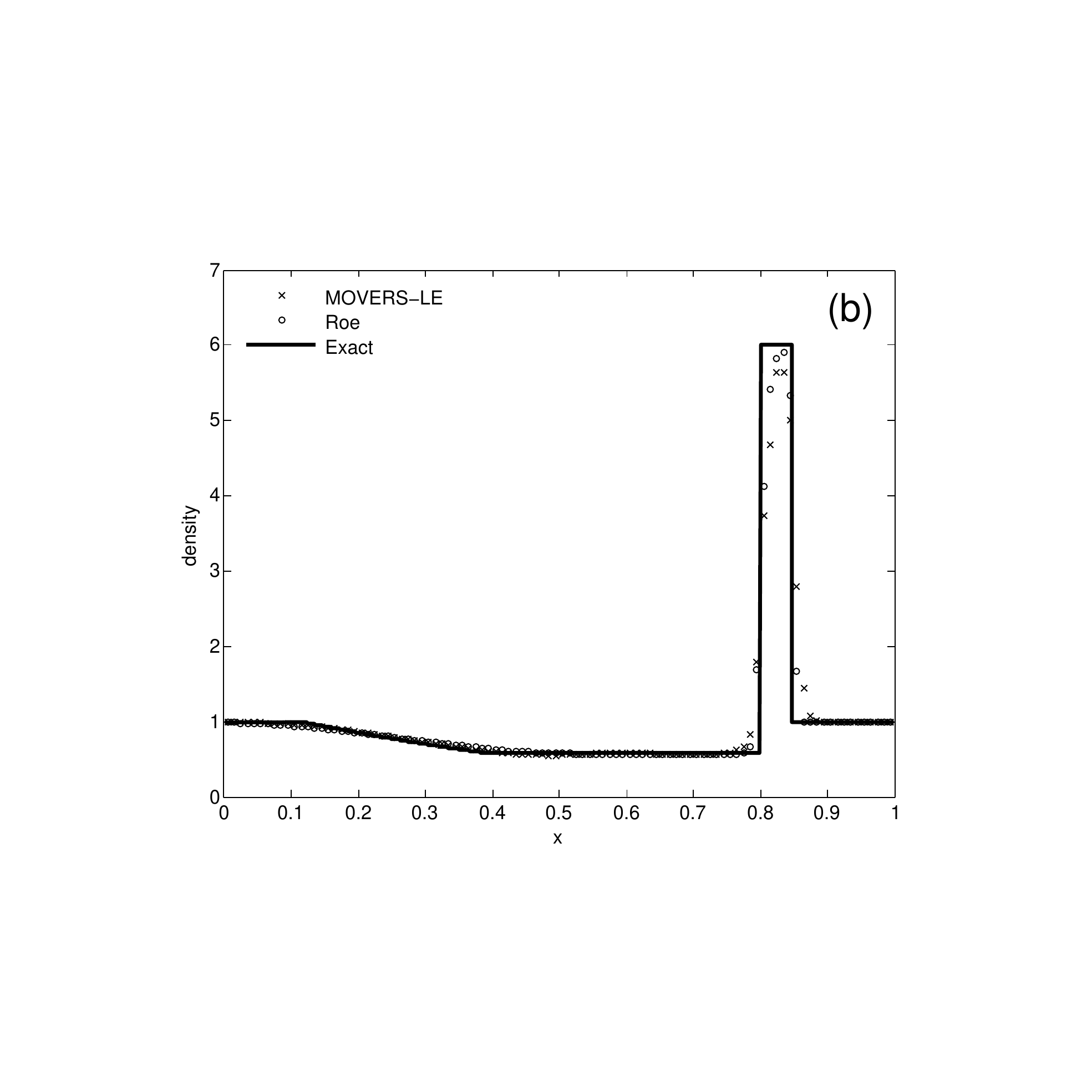}}
 \caption{(a) MOVERS-LE for Strong Discontinuities (b) Comparison of MOVERS-LE and Roe scheme for Slowly Moving Contact Discontinuity}
 \label{FigMoversLEProb45}
\end{figure}

This is severe test case as the solution of the problem has three discontinuities. This test case is also used to test the numerical method for accuracy and robustness. As shown in the Fig(\ref{FigMoversLEProb45}a), MOVERS-LE does not encounter any problem and is robust.  

\subsection{Slowly moving contact}
This test case is designed to test the numerical algorithm for its ability to resolve slowly moving contact discontinuities.   The solution of this test problem has a shock wave and a slowly moving contact discontinuity to right and a rarefaction wave to the left. As shown in the Fig(\ref{FigMoversLEProb45}b), the MOVERS-LE captures the slowly moving contact discontinuity quite accurately and its accuracy is comparable to that or Roe scheme.   

\subsection{Test cases for 2D-Euler equations}
The present scheme, MOVERS-LE, is tested on various benchmark problems for 2D-Euler equations. The problems are chosen to assess the numerical scheme for accuracy, robustness and ability to avoid shock instabilities.  

\subsubsection{Grid aligned slip flow}
\begin{figure}[h!]
 \subfigure{\includegraphics[trim=80.0 30.0 80.0 55.0, clip, width=0.5\textwidth]{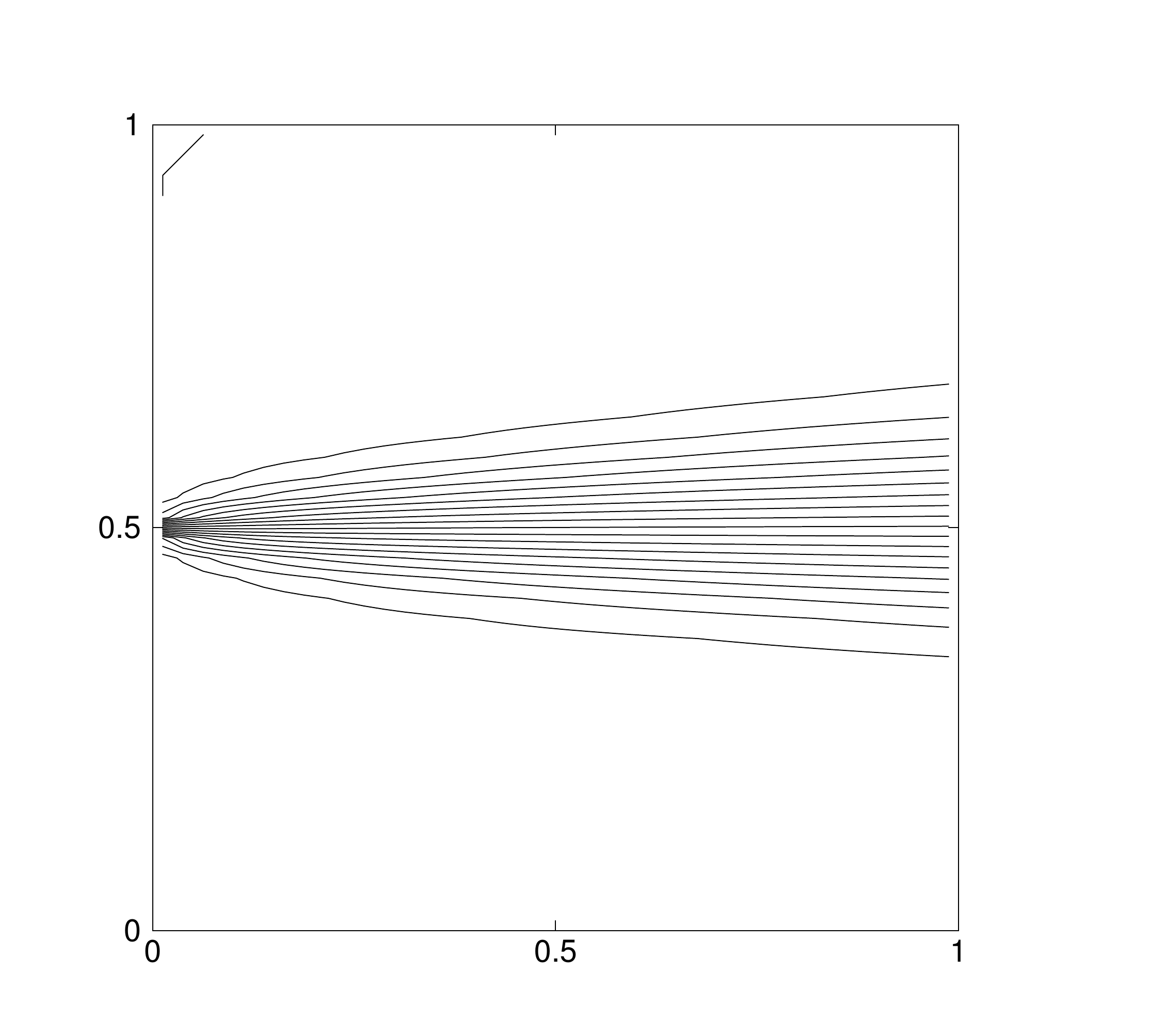}}
 \subfigure{\includegraphics[trim=80.0 30.0 80.0 55.0, clip, width=0.5\textwidth]{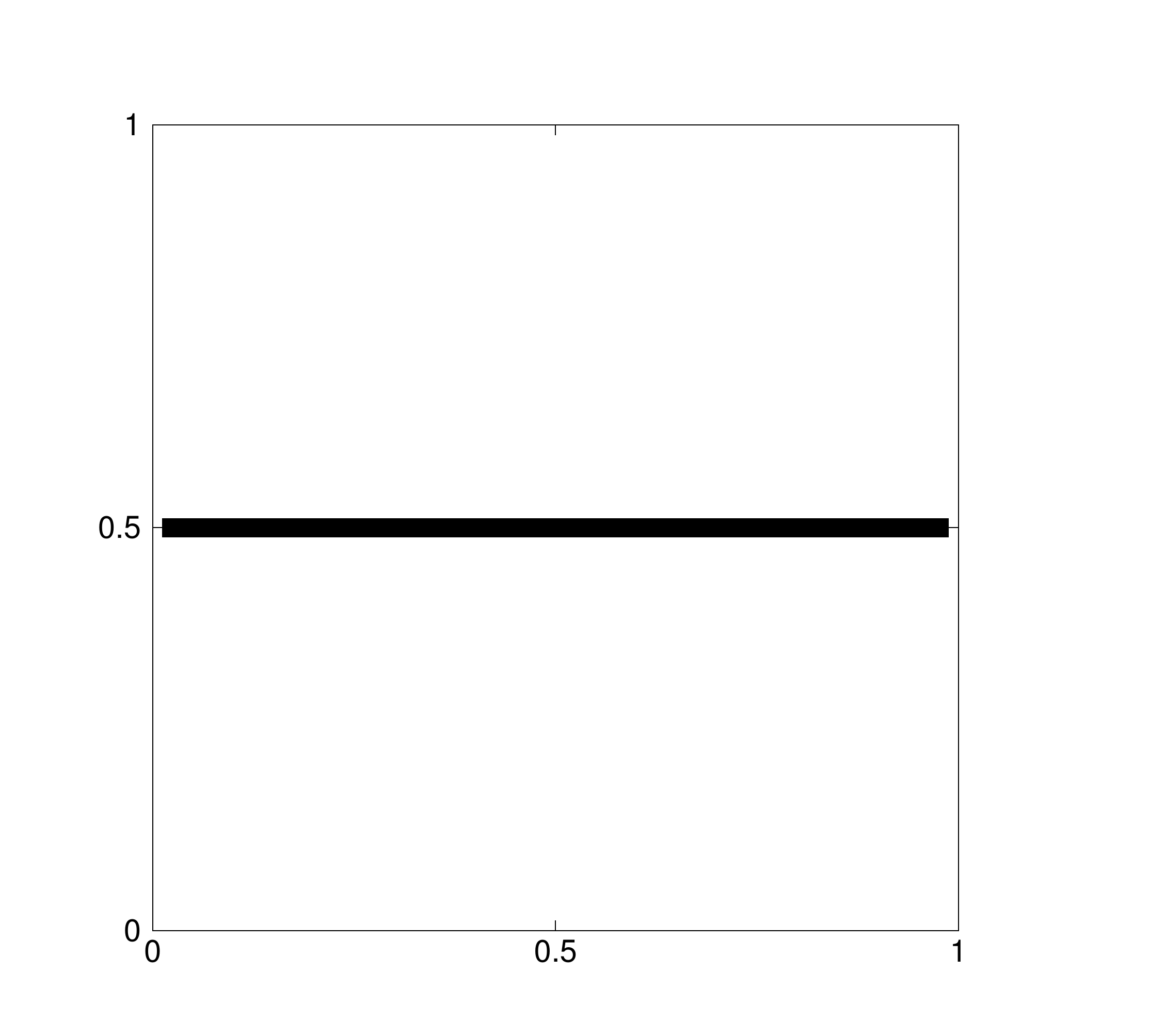}}
 \caption{Comparison of grid aligned slip flow problem results for MOVERS-LE (right) with LLF method (left) on a grid with 40 $\times$ 40 grid points: Mach contours (2.0: 0.05 :3.0)}
 \label{FigSlipFlow}
\end{figure}

In this test case\cite{manna1992three}, a Mach $3$ flow slips over a Mach $2$ flow, with the pressure gradient being zero and no jump in density.  This slip flow problem tests the schemes for its ability to exactly capture grid aligned discontinuities. The left side of the domain is prescribed with inflow values and  at all other boundaries, values are extrapolated from inside. Fig. (\ref{FigSlipFlow}) shows the results with the present scheme, which exactly captures the slip stream as compared with LLF method. 

\subsubsection{Oblique shock reflection}
\begin{figure}[h!]
 \subfigure{\includegraphics[trim=2.0 45.0 50.0 345.0, clip, width=0.5\textwidth]{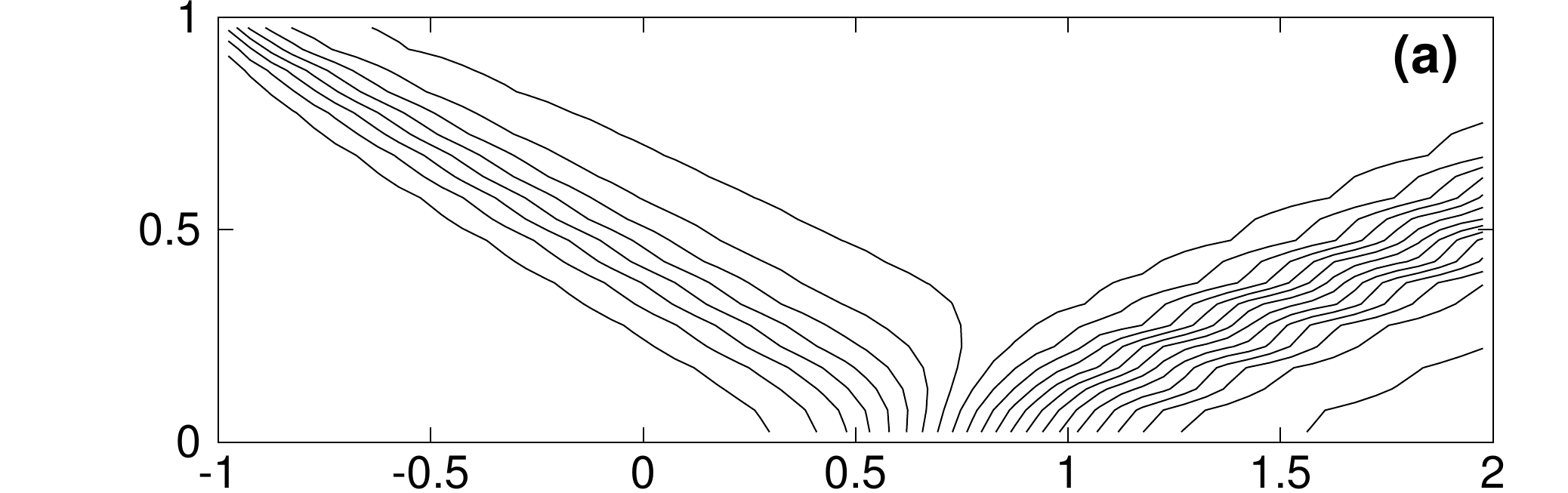}}
 \subfigure{\includegraphics[trim=2.0 45.0 50.0 345.0, clip, width=0.5\textwidth]{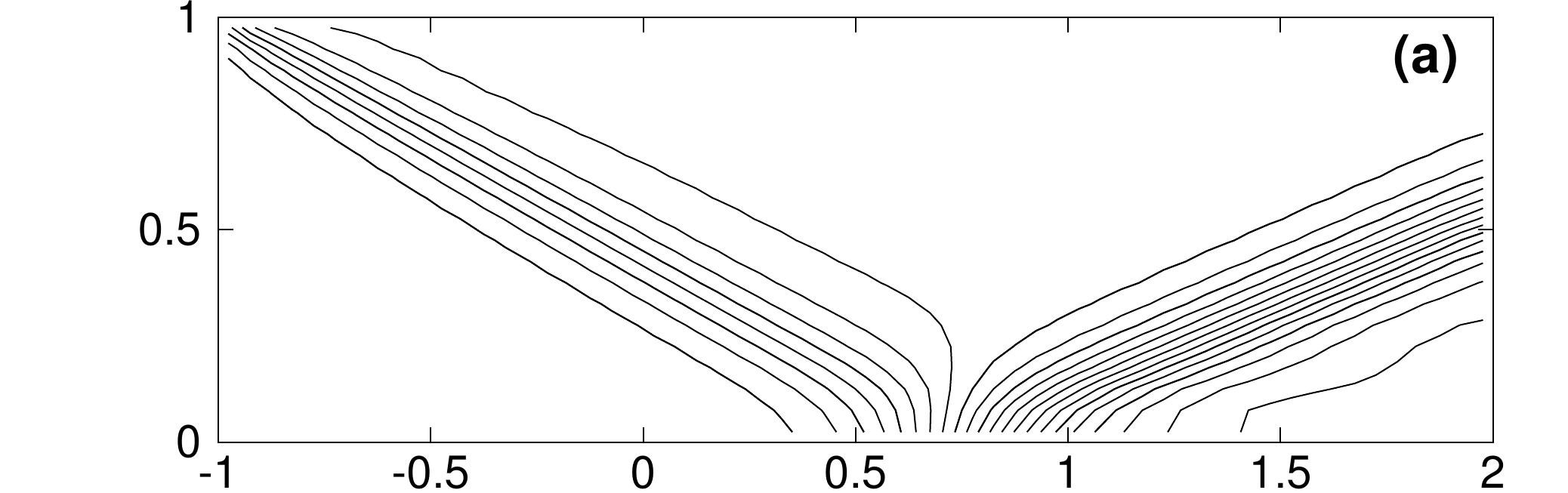}}
 \subfigure{\includegraphics[trim=2.0 45.0 50.0 345.0, clip, width=0.5\textwidth]{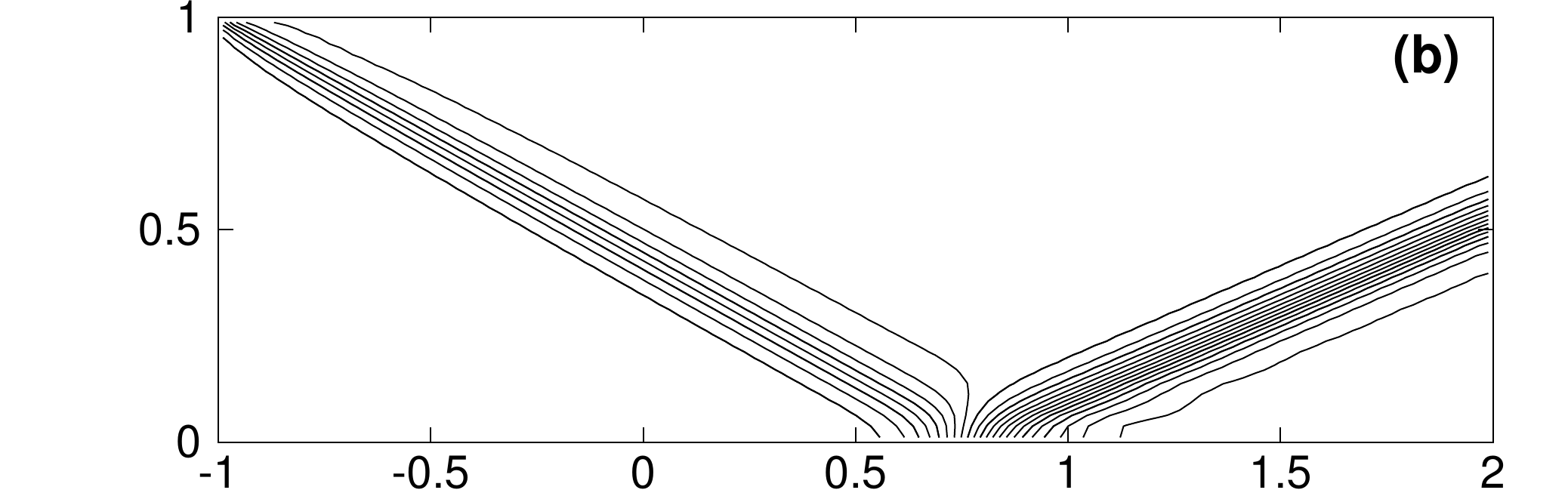}}
 \subfigure{\includegraphics[trim=2.0 45.0 50.0 345.0, clip, width=0.5\textwidth]{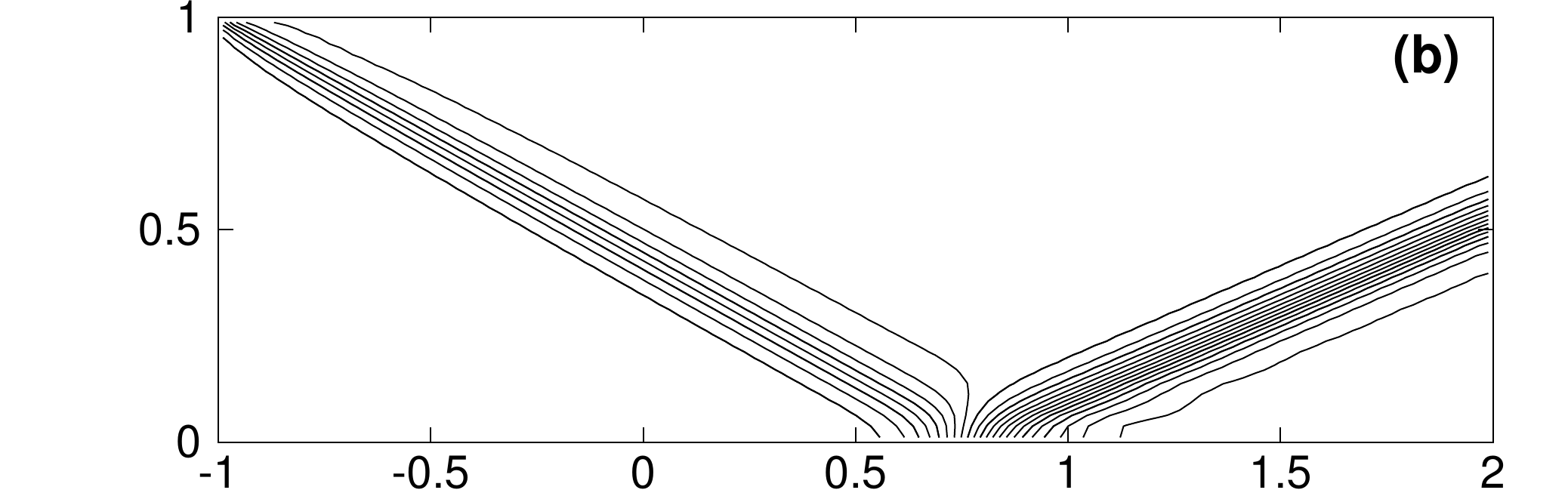}}
 \subfigure{\includegraphics[trim=2.0 45.0 50.0 345.0, clip, width=0.5\textwidth]{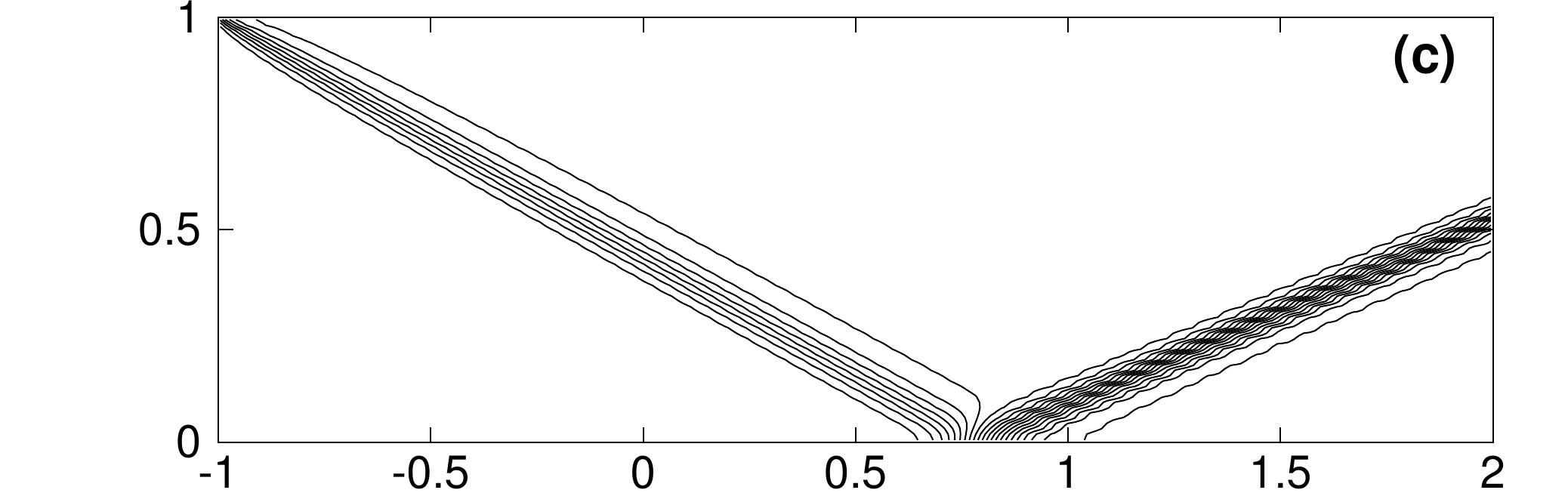}}
 \subfigure{\includegraphics[trim=2.0 45.0 50.0 345.0, clip, width=0.5\textwidth]{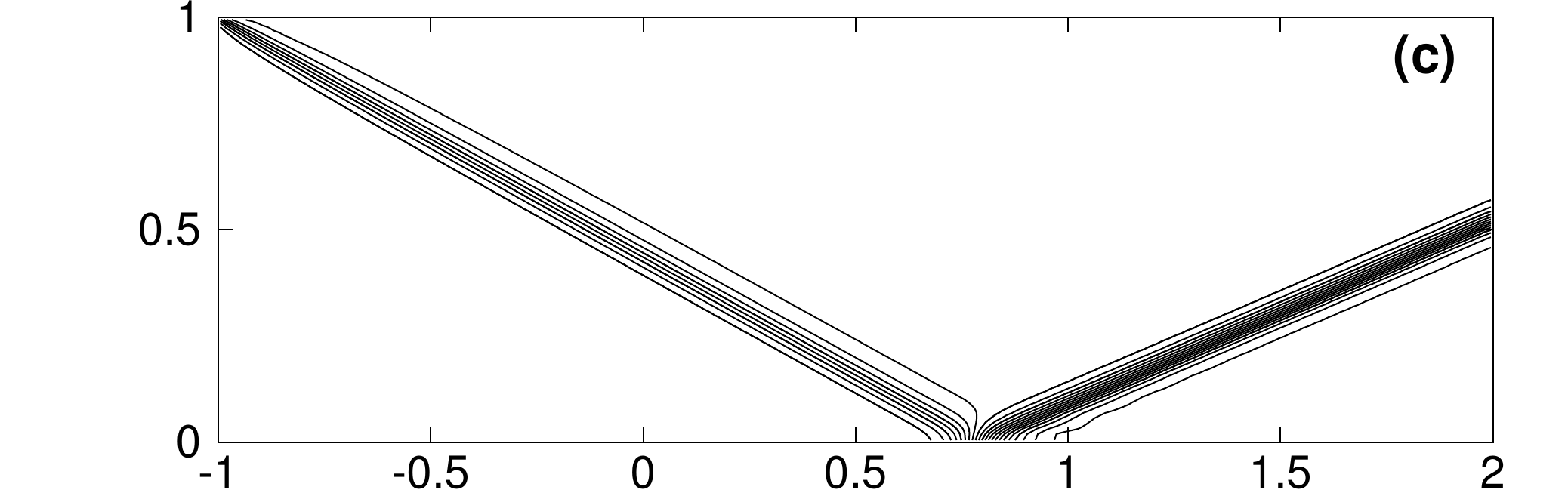}}
 \caption{Comparison of first order results of MOVERS (left) with MOVERS-LE (right): pressure contours (0.7:0.1:2.9) for the grids: (a) 60$\times$20, (b) 120$\times$40 and (c) 240$\times$80 }
 \label{FigShkReflFo}
\end{figure}

\begin{figure}[h!]
 \subfigure{\includegraphics[trim=2.0 45.0 50.0 345.0, clip, width=0.5\textwidth]{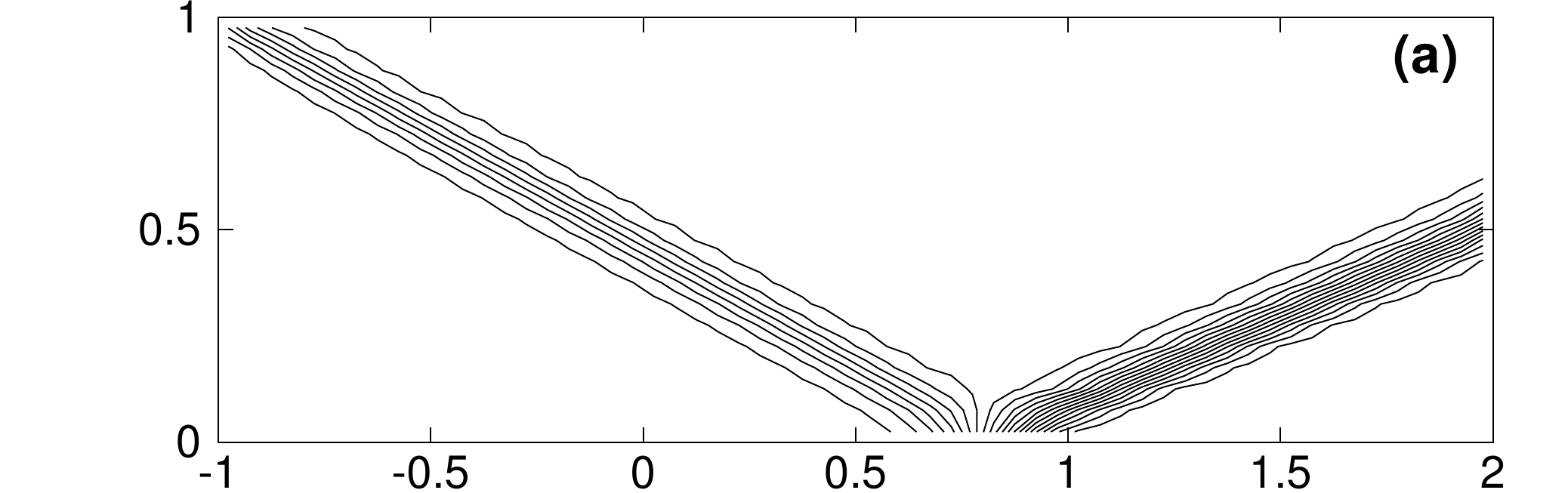}}
 \subfigure{\includegraphics[trim=2.0 45.0 50.0 345.0, clip, width=0.5\textwidth]{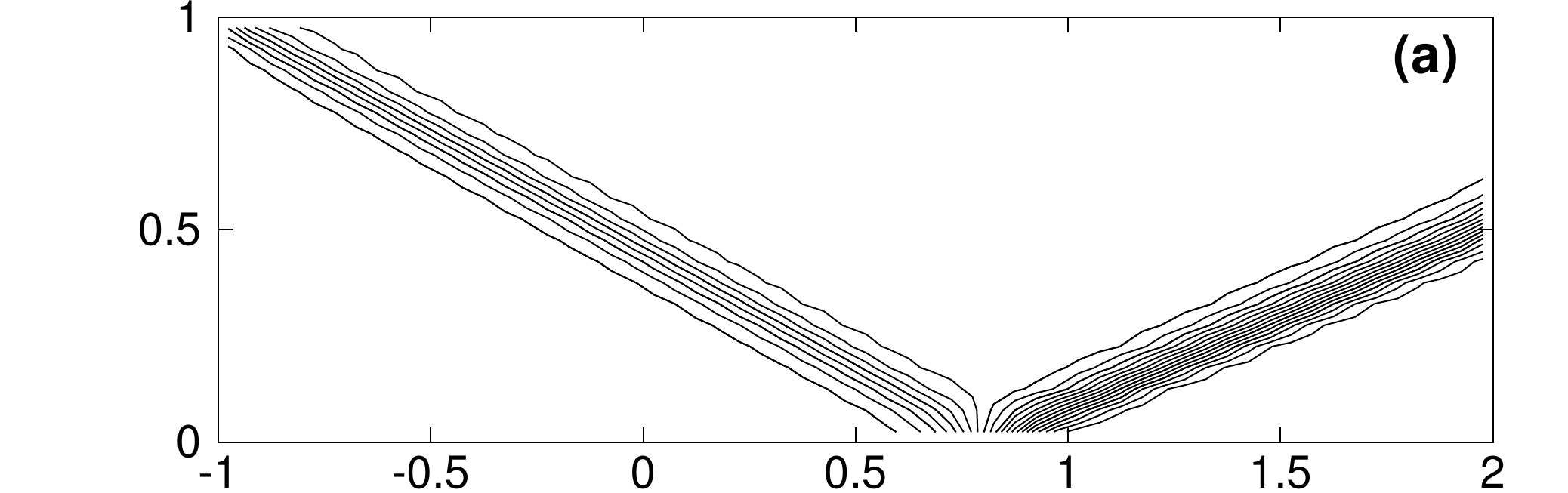}}
 \subfigure{\includegraphics[trim=2.0 45.0 50.0 345.0, clip, width=0.5\textwidth]{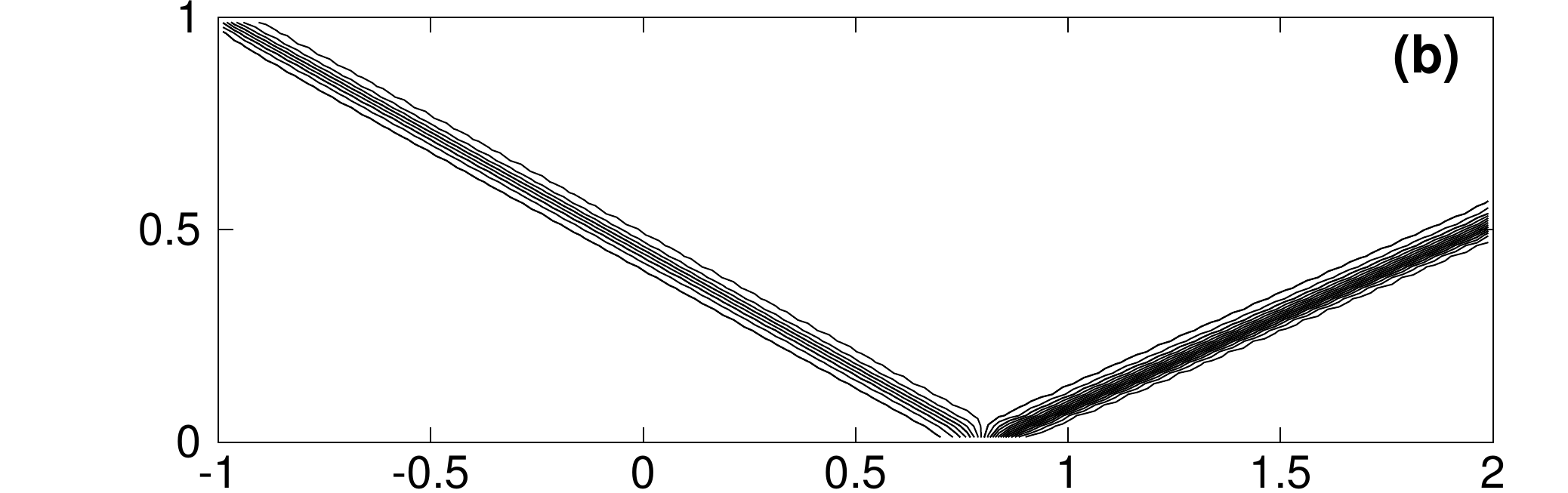}}
 \subfigure{\includegraphics[trim=2.0 45.0 50.0 345.0, clip, width=0.5\textwidth]{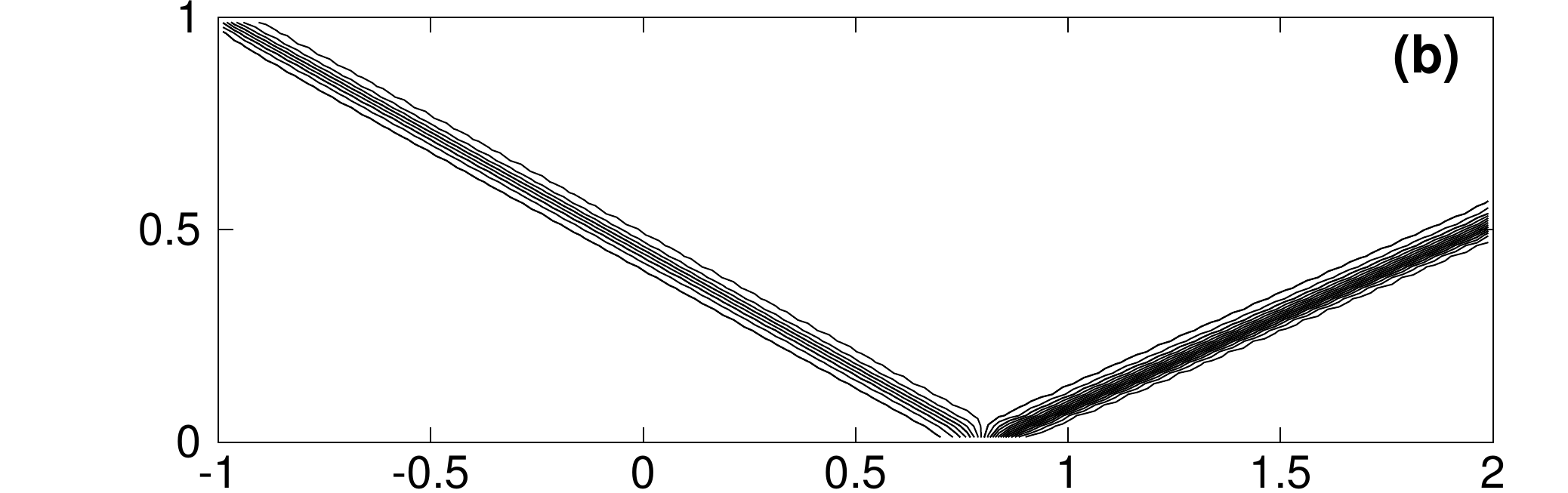}}
 \subfigure{\includegraphics[trim=2.0 45.0 50.0 345.0, clip, width=0.5\textwidth]{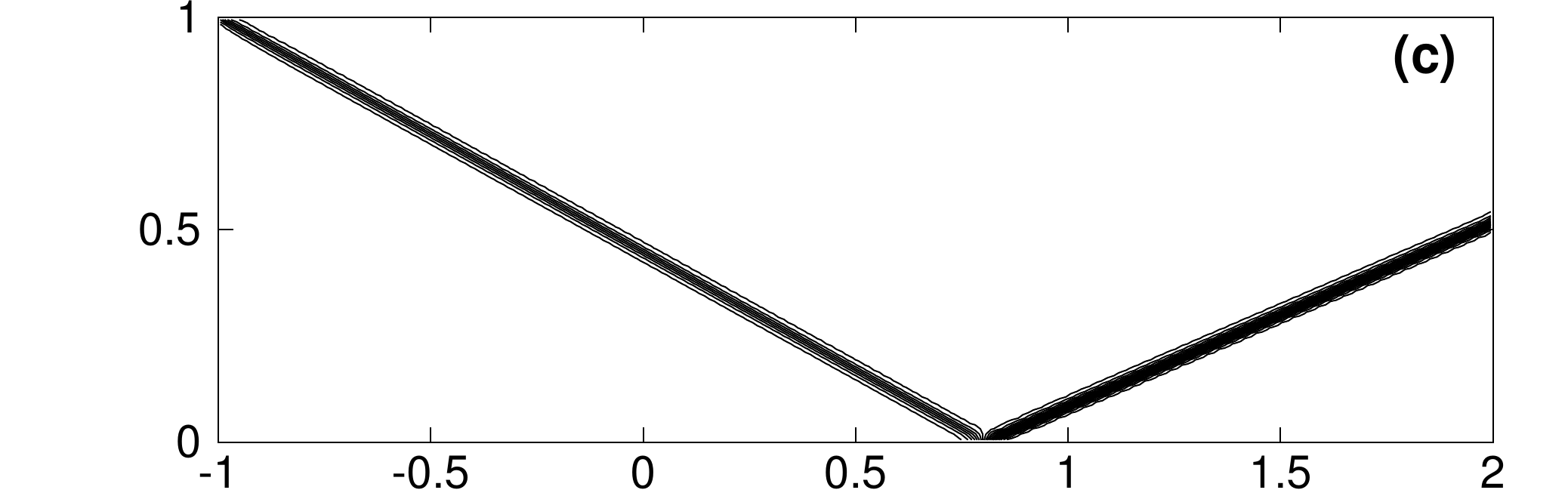}}
 \subfigure{\includegraphics[trim=2.0 45.0 50.0 345.0, clip, width=0.5\textwidth]{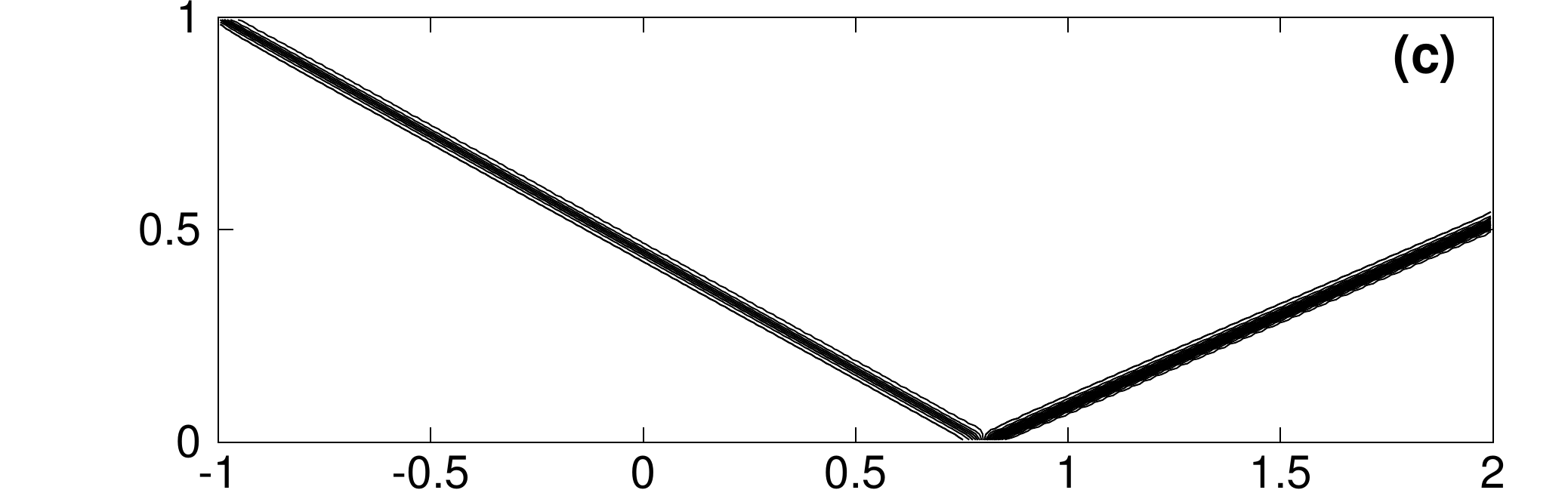}}
 \caption{Comparison of second order results of MOVERS (left) with MOVERS-LE (right): pressure contours (0.7:0.1:2.9) for the grids: (a) 60$\times$20, (b) 120$\times$40 and (c) 240$\times$80 }
 \label{FigShkReflSo}
\end{figure}

In this test case \cite{jin1995relaxation} an oblique shock is introduced at the top left corner of the computational domain. This oblique shock hits the solid wall and reflects back. The initial conditions for this problem are as given below.
\begin{equation*}
\left( \rho, u, v, p \right)_{0,y,t}=\left( 1.0, 2.9, 0, 1/1.4 \right) 
\end{equation*}
\begin{equation*}
\left( \rho, u, v, p \right)_{x,1,t}=\left( 1.69997, 2.61934, -0.50633, 1.52819 \right)
\end{equation*}
At the left and top sides of the domain inflow and post shock boundary conditions are prescribed as given above.  Solid-wall and supersonic outflow boundary conditions are prescribed at the bottom and right boundaries.  The results for the oblique shock reflection with MOVERS-LE for both first and second order accuracy are shown in the Fig. (\ref{FigShkReflFo}) and Fig. (\ref{FigShkReflSo}) respectively.  As can be seen from the figures, the present scheme is slightly more accurate compared to MOVERS, and the contours are smoother.  


\subsubsection{Ramp in a channel} 
\begin{figure}[h!]
 \subfigure{\includegraphics[trim=2.0 45.0 50.0 345.0, clip, width=0.5\textwidth]{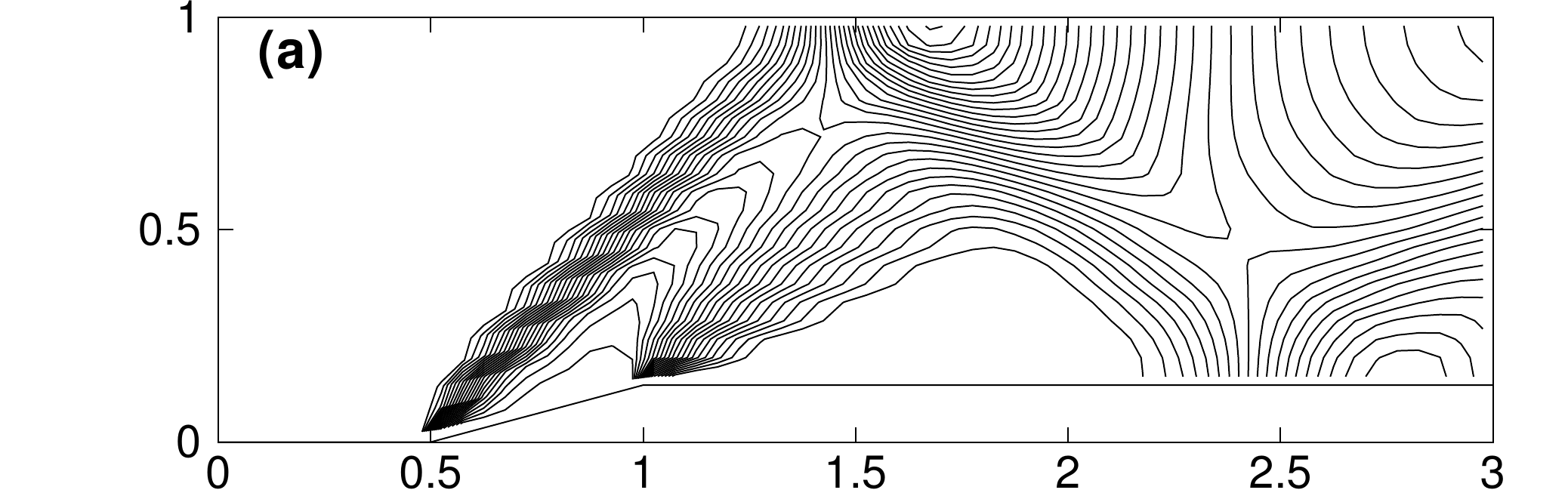}}
 \subfigure{\includegraphics[trim=2.0 45.0 50.0 345.0, clip, width=0.5\textwidth]{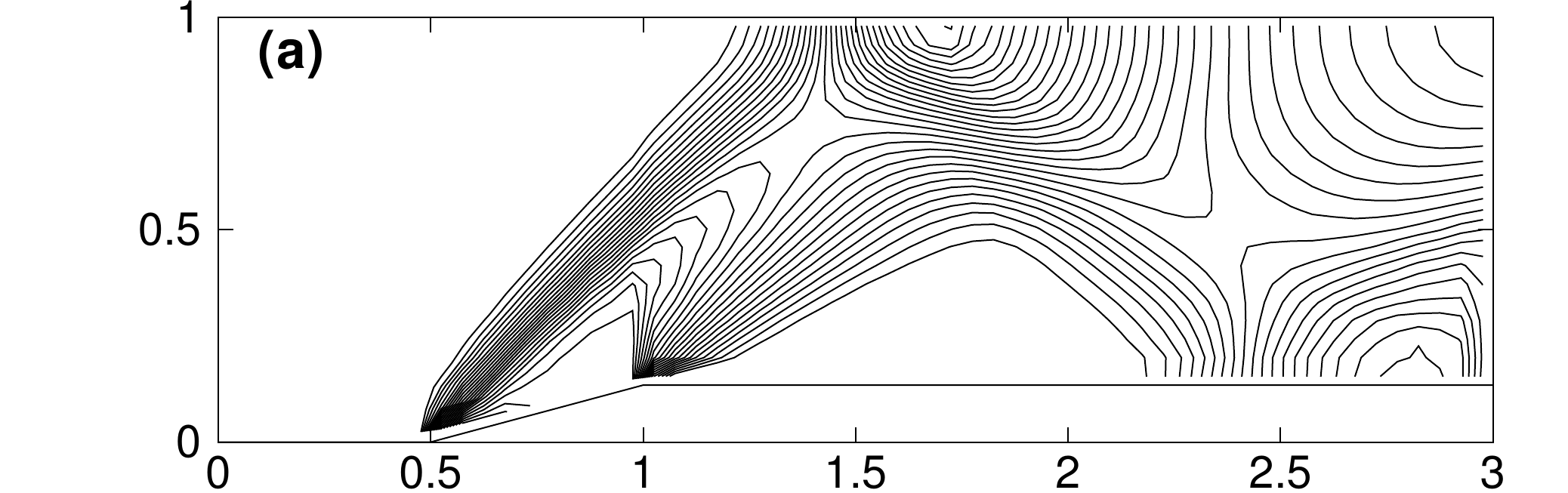}}
 \subfigure{\includegraphics[trim=2.0 45.0 50.0 345.0, clip, width=0.5\textwidth]{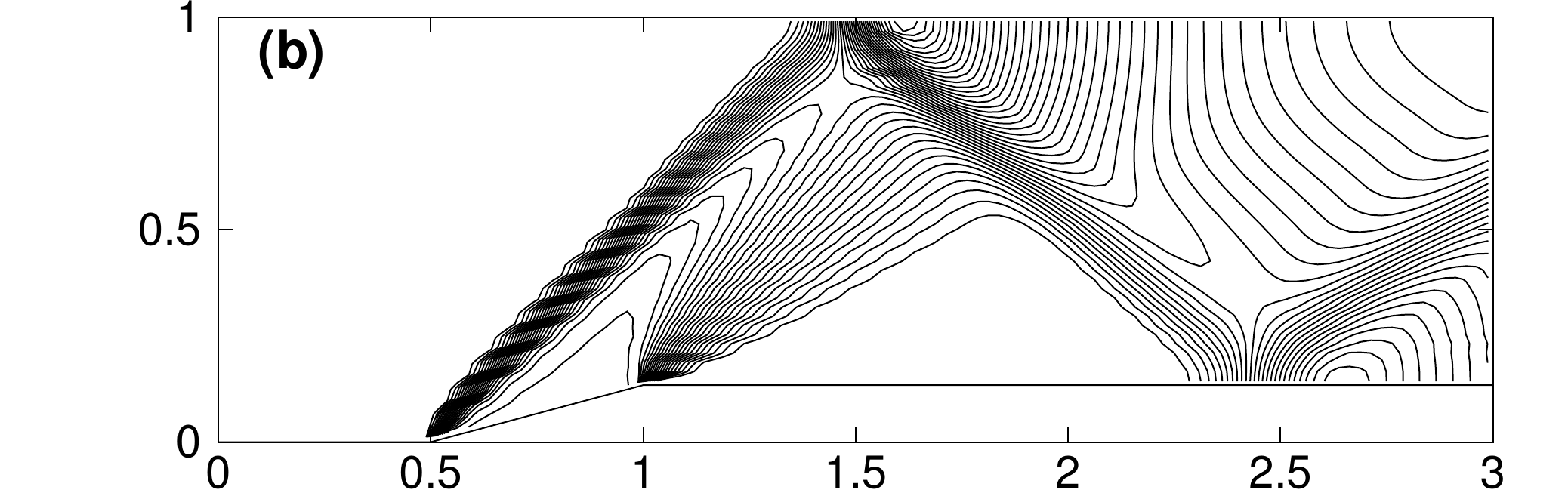}}
 \subfigure{\includegraphics[trim=2.0 45.0 50.0 345.0, clip, width=0.5\textwidth]{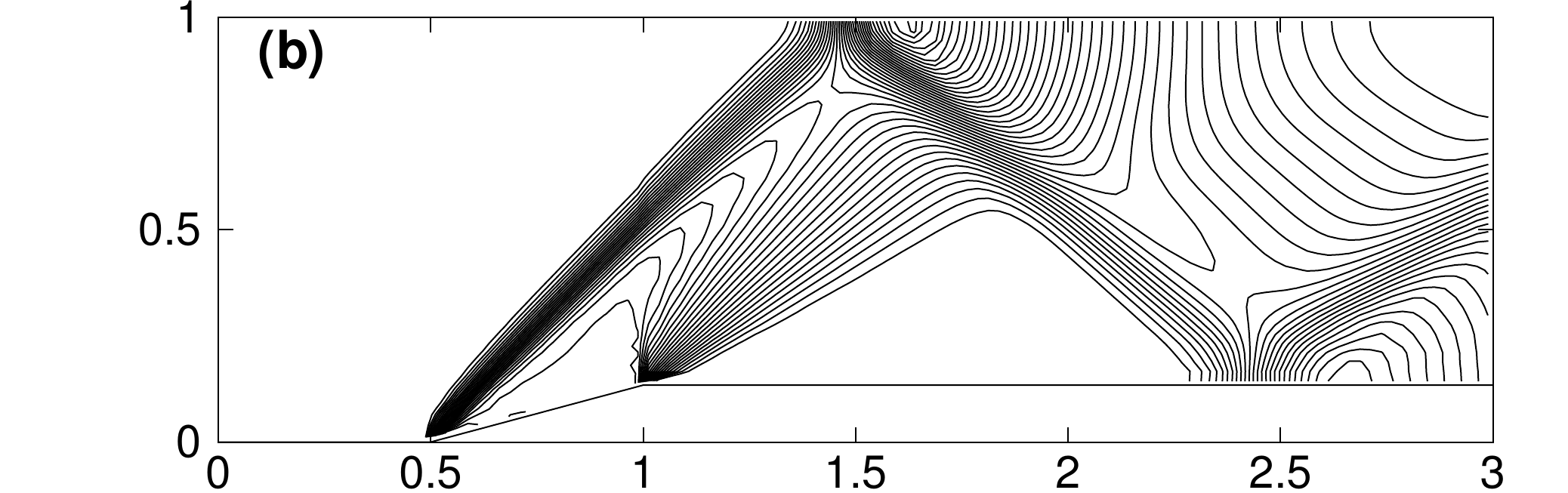}}
 \subfigure{\includegraphics[trim=2.0 45.0 50.0 345.0, clip, width=0.5\textwidth]{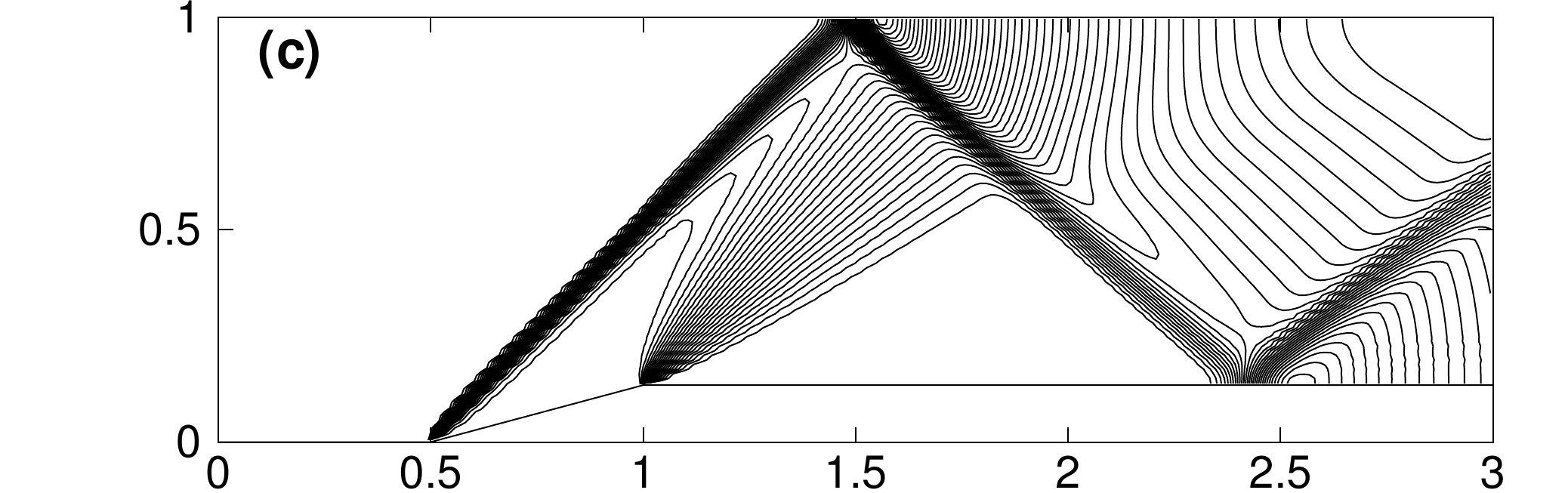}}
 \subfigure{\includegraphics[trim=2.0 45.0 50.0 345.0, clip, width=0.5\textwidth]{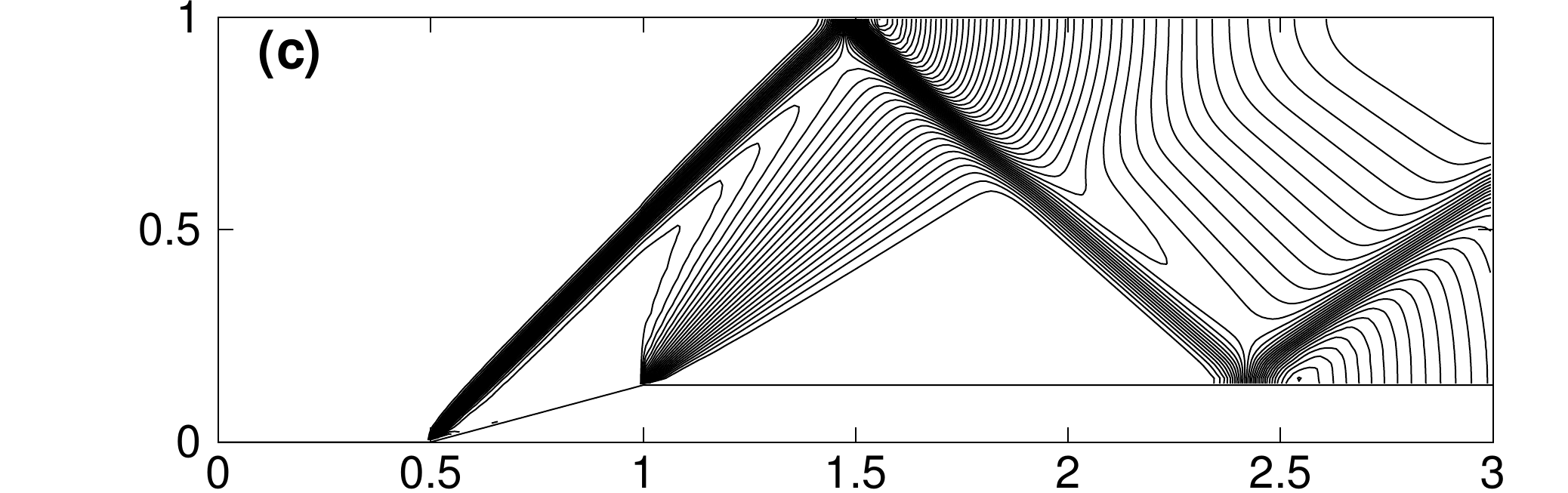}}
 \caption{Comparison of first order results of MOVERS (left) with MOVERS-LE (right): pressure contours (0.7:0.1:2.9) for the grids: (a) 60$\times$20, (b) 120$\times$40 and (c) 240$\times$80 }
 \label{FigRampFo}
\end{figure}

\begin{figure}[h!]
 \subfigure{\includegraphics[trim=2.0 45.0 50.0 345.0, clip, width=0.5\textwidth]{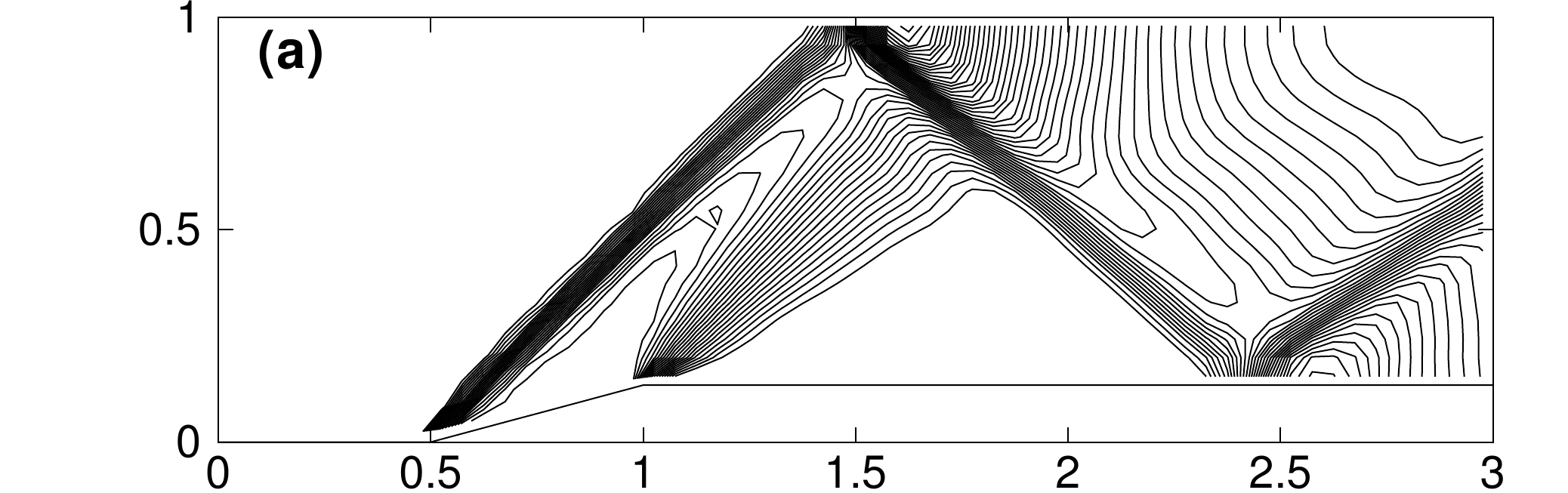}}
 \subfigure{\includegraphics[trim=2.0 45.0 50.0 345.0, clip, width=0.5\textwidth]{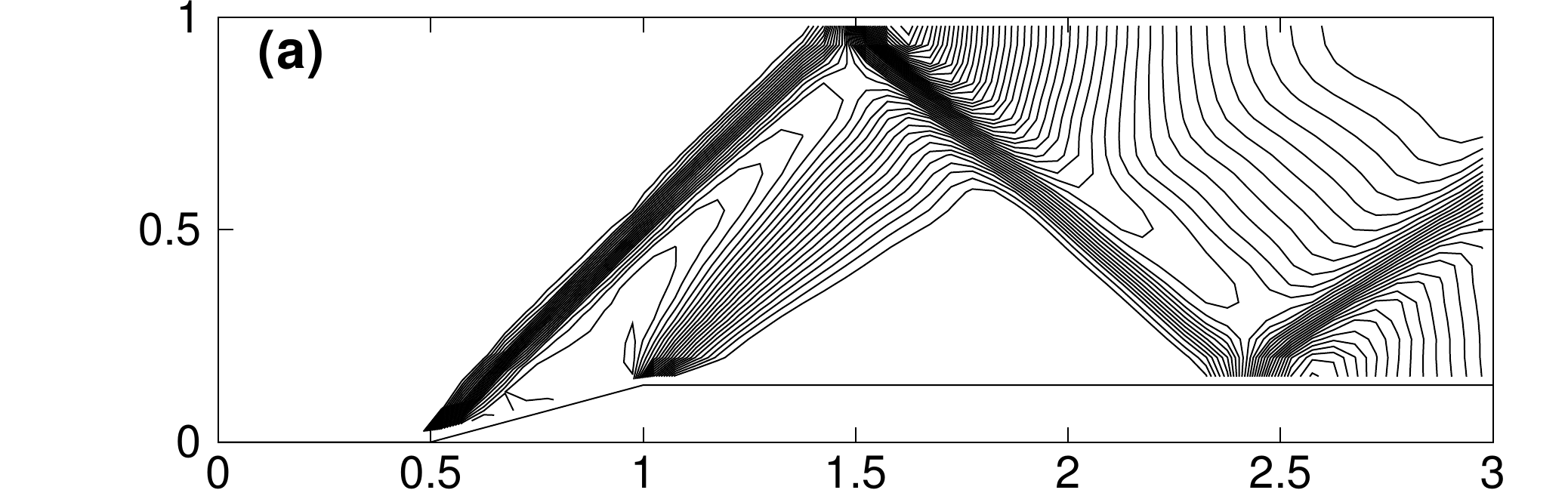}}
 \subfigure{\includegraphics[trim=2.0 45.0 50.0 345.0, clip, width=0.5\textwidth]{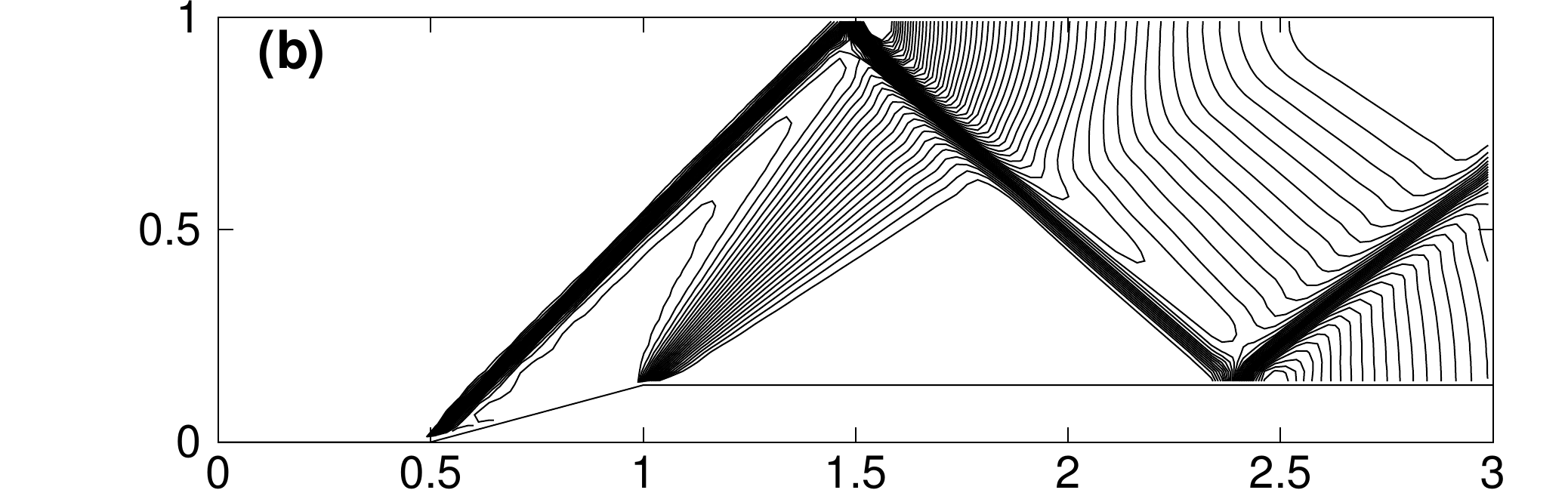}}
 \subfigure{\includegraphics[trim=2.0 45.0 50.0 345.0, clip, width=0.5\textwidth]{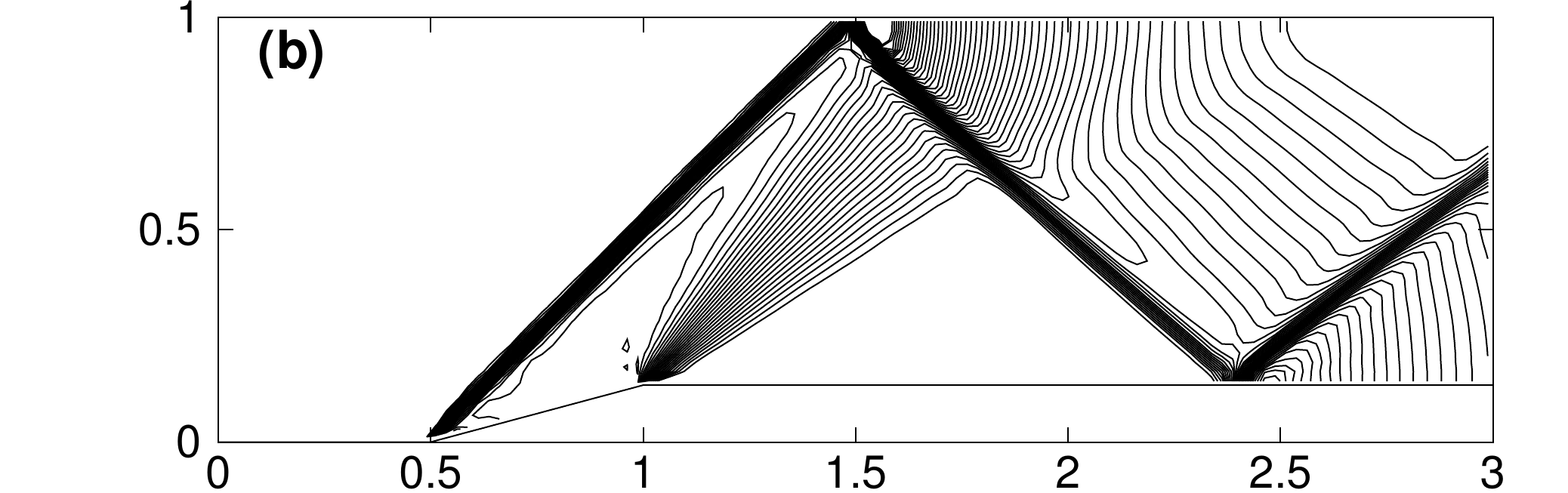}}
 \subfigure{\includegraphics[trim=2.0 45.0 50.0 345.0, clip, width=0.5\textwidth]{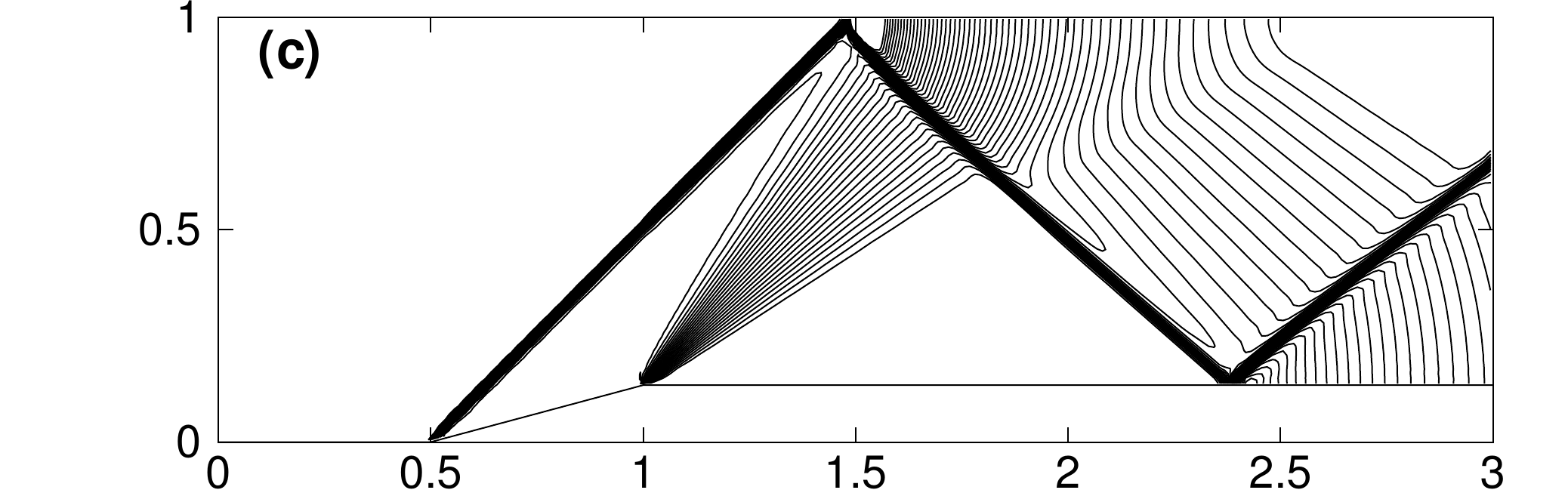}}
 \subfigure{\includegraphics[trim=2.0 45.0 50.0 345.0, clip, width=0.5\textwidth]{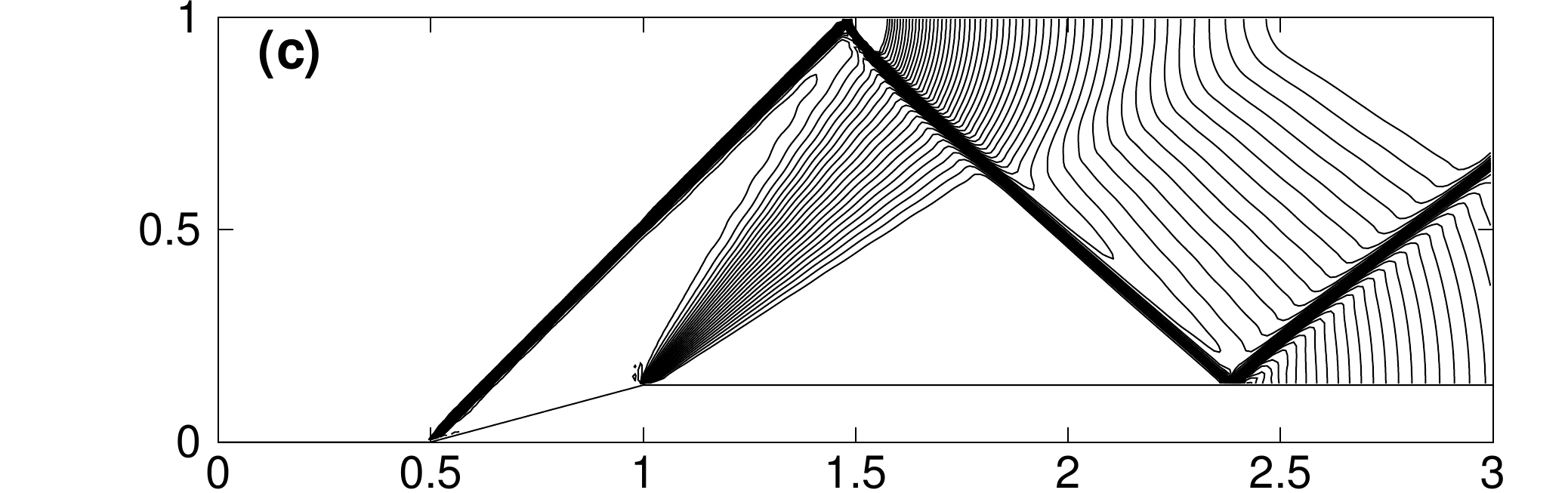}}
 \caption{Comparison of second order results of MOVERS (left) with MOVERS-LE (right): pressure contours (0.7:0.1:2.9) for the grids: (a) 60$\times$20, (b) 120$\times$40 and (c) 240$\times$80 }
 \label{FigRampSo}
\end{figure}

In this test case~\citep{levy1993use} a Mach $2$ flow encounters a wind tunnel with a $15^0$ ramp.  The flow feature of this problem includes an oblique shock which reflects from the upper wall and interacts with the the expansion fan generated at the corner of the ramp.  Fig. (\ref{FigRampFo}) and Fig. (\ref{FigRampSo}) show the first and second order accurate results of MOVERS-LE compared with MOVERS.  Although there is no noticeable difference in accuracy, it is worth noting that the contours are smoother with the use of present scheme.  


\subsubsection{Hypersonic flow past a half cylinder}
\begin{figure}[h!]
 \centering
 \textbf{}
 \subfigure{\includegraphics[trim=180.0 100.0 300.0 100.0, clip, width=0.2\textwidth]{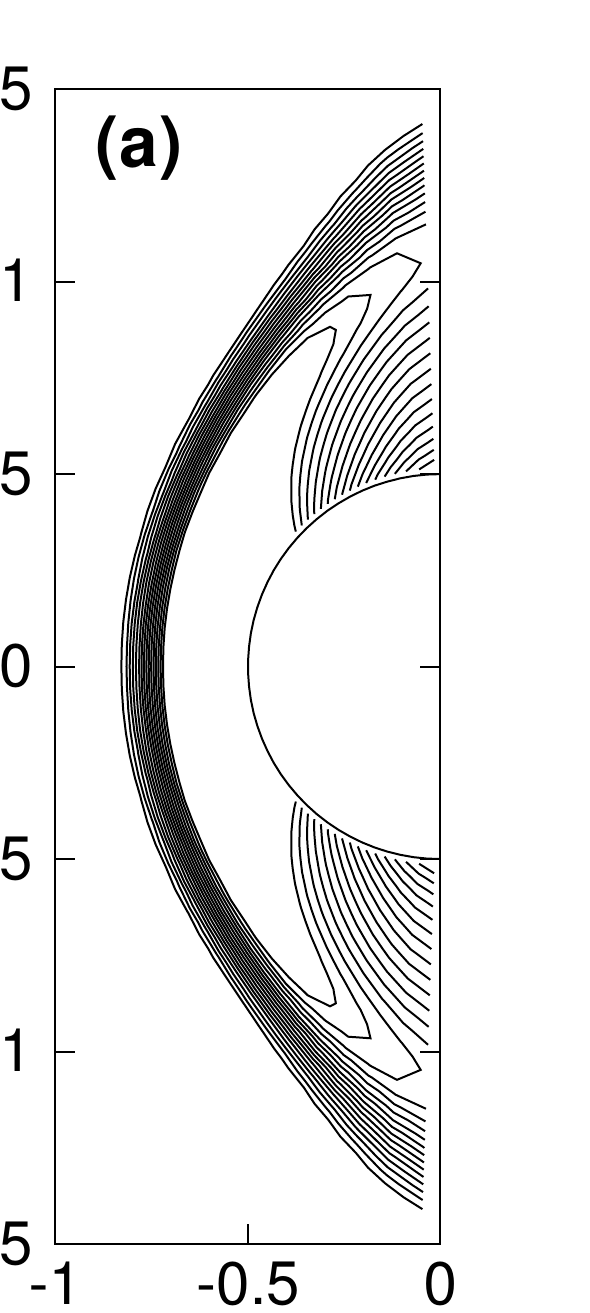}}
 \subfigure{\includegraphics[trim=180.0 100.0 300.0 100.0, clip, width=0.2\textwidth]{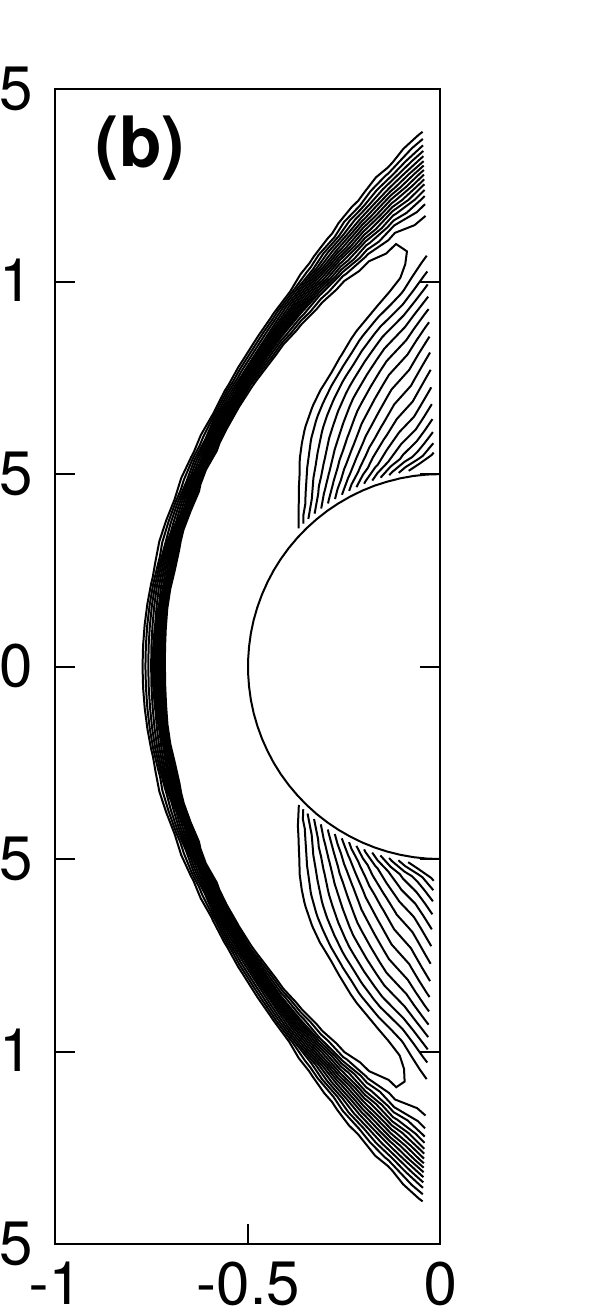}}
 \subfigure{\includegraphics[trim=180.0 100.0 300.0 100.0, clip, width=0.2\textwidth]{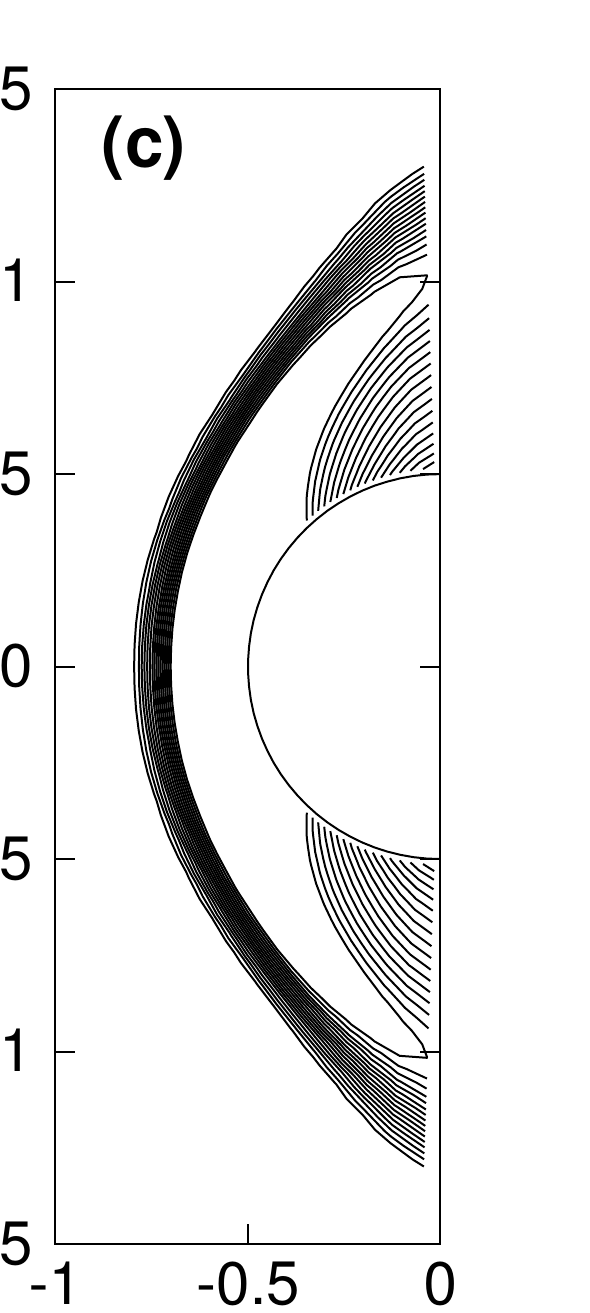}}
 \subfigure{\includegraphics[trim=180.0 100.0 300.0 100.0, clip, width=0.2\textwidth]{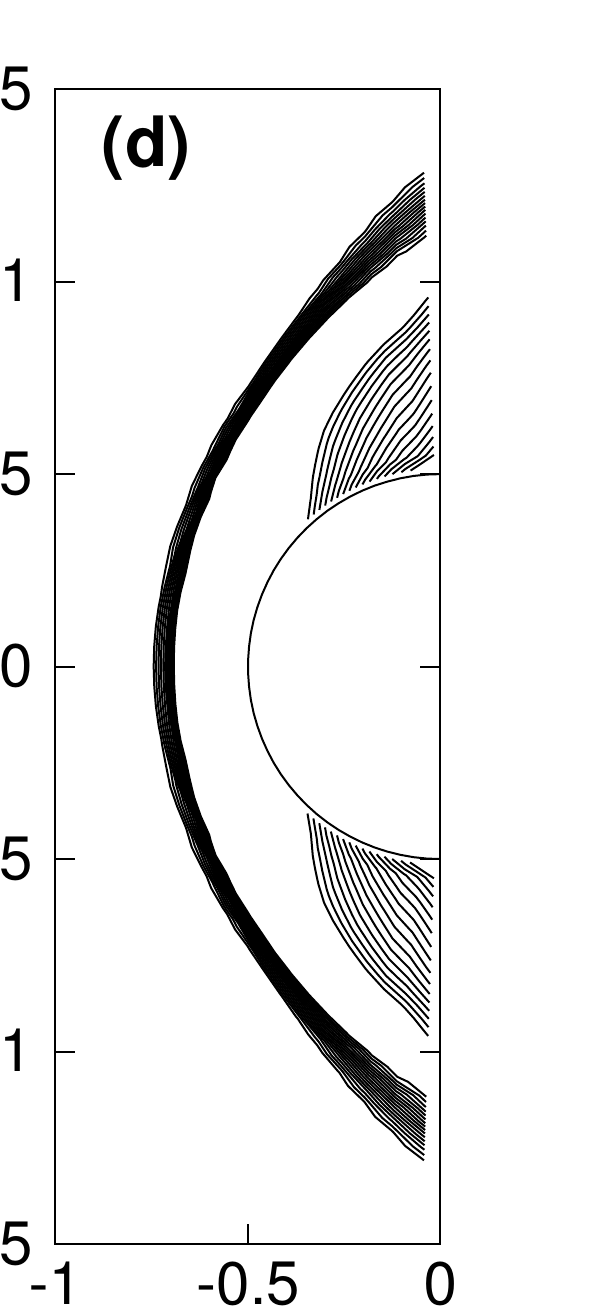}}
 \caption{Comparison of First order results of MOVERS-LE with LLF on 45 $\times 45$ grids: density contours (2.0: 0.2: 5.0) (a) LLF for Mach $6$ flow (b) MOVERS-LE for Mach $6$ flow (c) LLF for Mach $20$ flow (d) MOVERS-LE Mach $20$ flow}
 \label{FigCylinFo}
\end{figure}

This is one of the multidimensional problems used to test the numerical algorithm for shock instability.  Many upwind schemes generate the infamous carbuncle shock as demonstrated by Quirk \cite{quirk1994a} and Meng-Sing Liou \cite{liou2000mass}.  As shown in the Fig. (\ref{FigCylinFo}), the present method does not encounter any problem in capturing the flow features.  

\subsubsection{Forward facing step}
\begin{figure}[h!]
 \centering
 \textbf{}
 \subfigure{\includegraphics[trim=80.0 70.0 135.0 350.0, clip, width=0.49\textwidth]{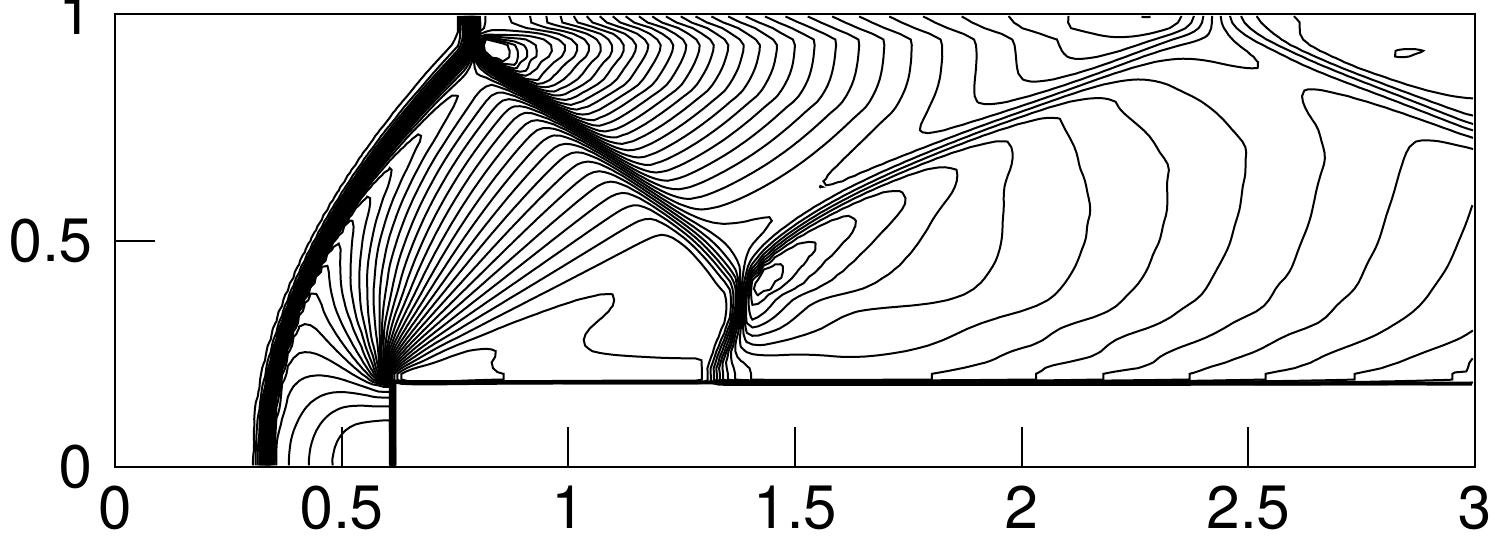}}
 \subfigure{\includegraphics[trim=80.0 70.0 135.0 350.0, clip, width=0.49\textwidth]{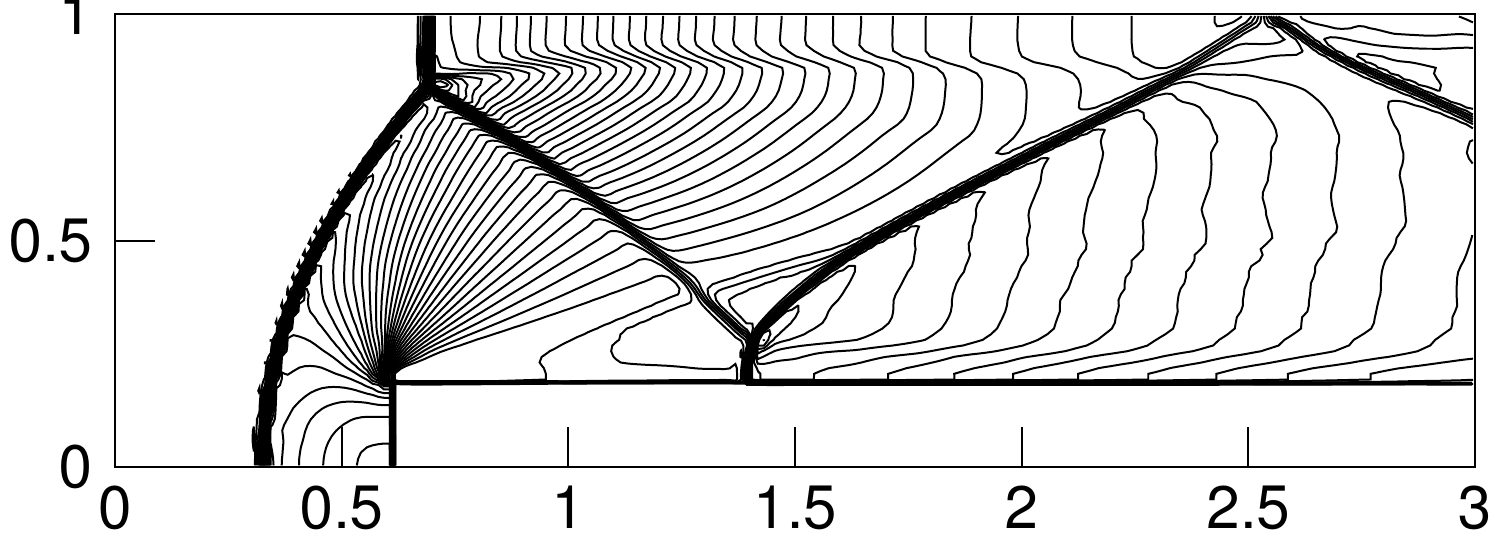}}
 \caption{Density contours (1.0: 6.5: 0.15) for forward facing step problem on 240 $\times$ 80 grids with MOVERS-LE for both first order (left) and second order (right) accuracy}
 \label{FigFFS}
\end{figure}

In this unsteady test case~\cite{woodward1984numerical}, a Mach $3$ flow enters the wind tunnel with a forward-facing step.   At time $t=4.0$ units, a lambda shock is visible with a triple point from which a slip surface is seen.   The reflected oblique shock from top wall interacts with the strong expansion generated at the corner of the step.  The left side of the computational domain is prescribed with inflow conditions.  Supersonic outflow boundary condition is enforced at the exit and wall boundary conditions are prescribed at top, bottom and on the step. Fig. (\ref{FigFFS}) shows the first and second order accurate results for this problem with MOVERS-LE.  As shown in the figure, the present scheme captures slip stream even with the first order accuracy, whereas many diffusive scheme like LLF do not capture it with first order accuracy.

\subsubsection{Shock diffraction}
\begin{figure}[h!]
 \centering
 \textbf{}
 \subfigure{\includegraphics[trim=70.0 46.0 120.0 50.0, clip, width=0.49\textwidth]{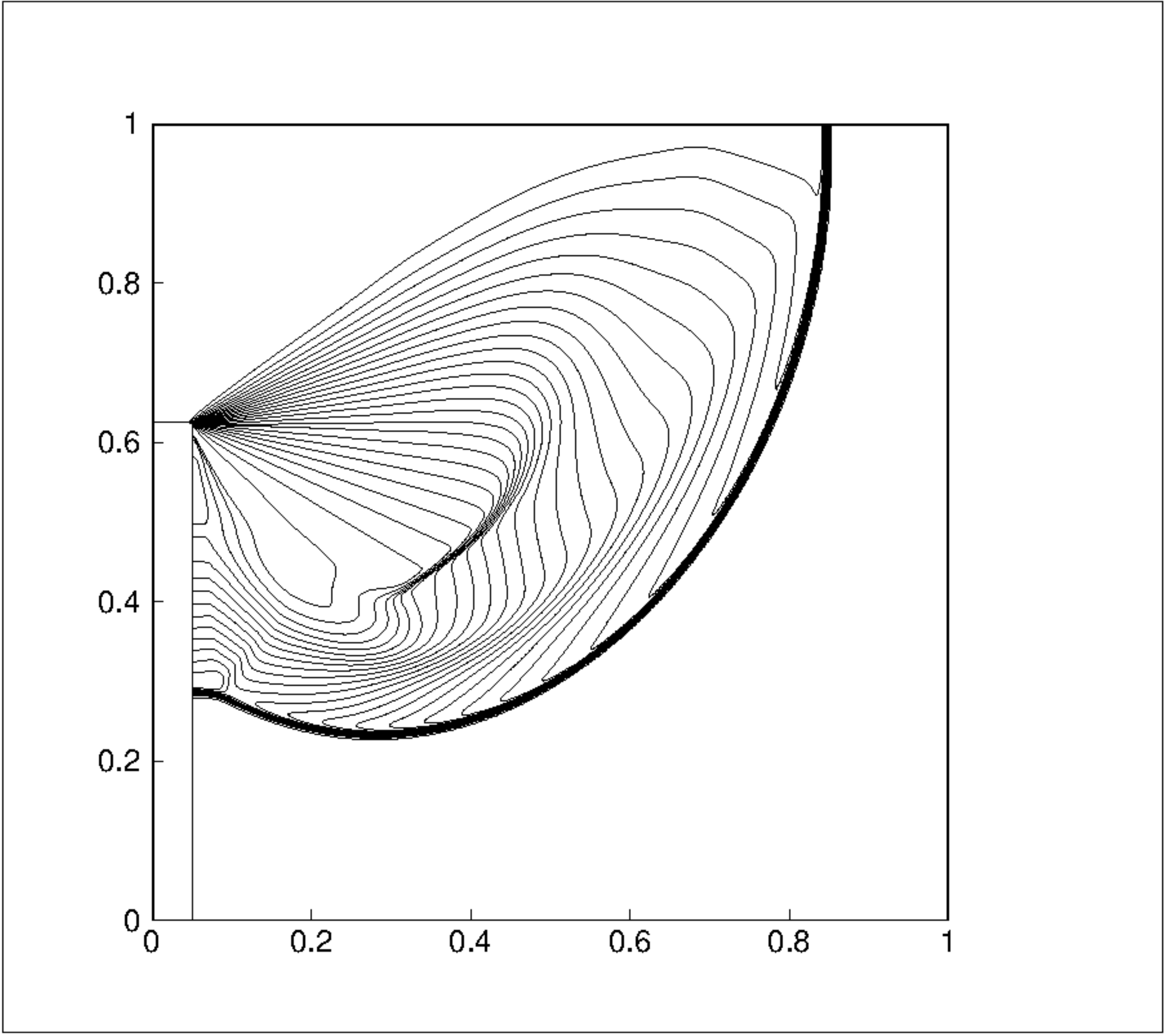}}
 \subfigure{\includegraphics[trim=70.0 46.0 120.0 50.0, clip, width=0.49\textwidth]{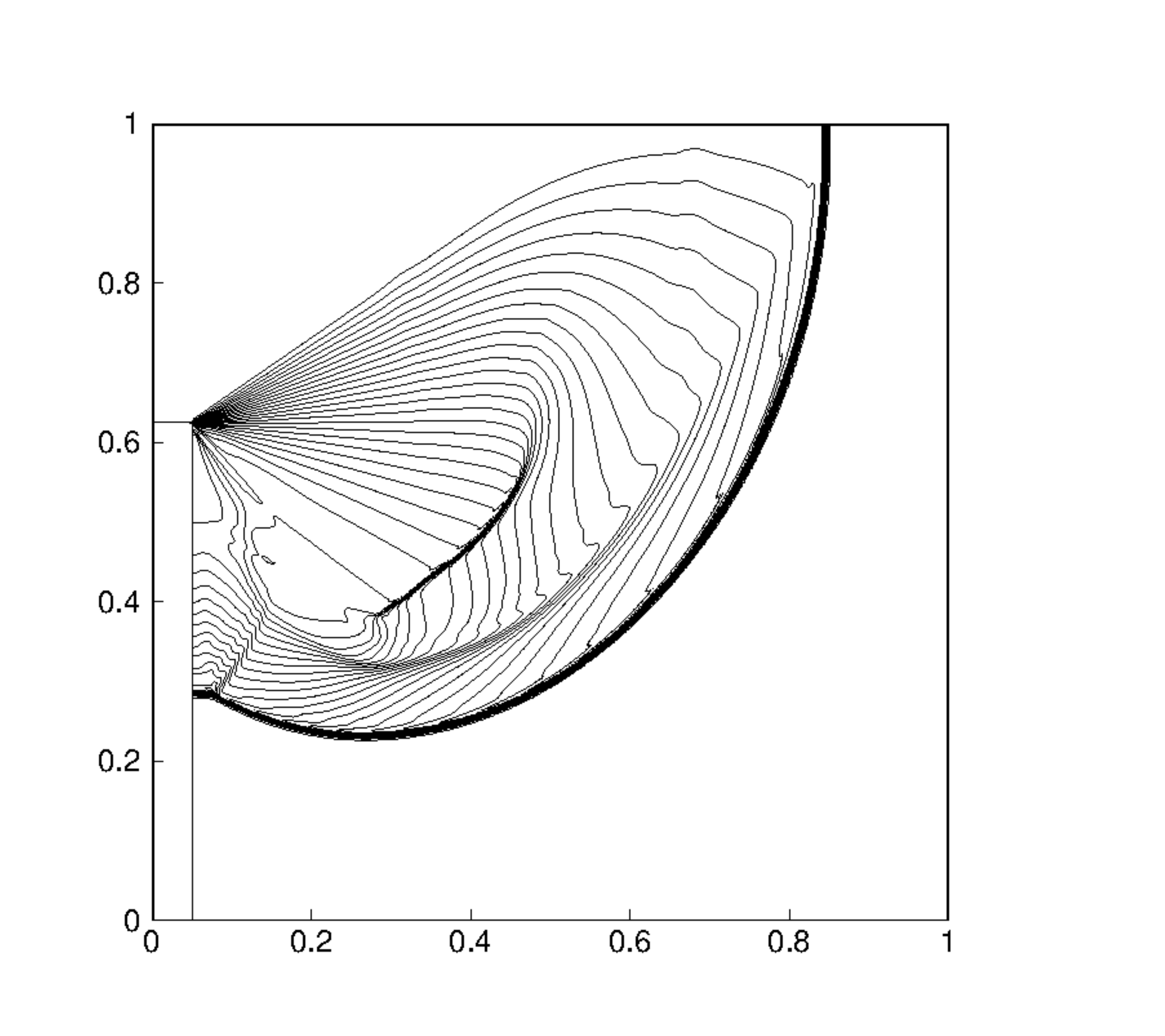}}
 \caption{Density contours (0.5: 7.0: 0.25) for shock diffraction problem on 400 $\times$ 400 grids with MOVERS-LE for both first order (left) and second order (right) accuracy}
 \label{FigShkDiffr}
\end{figure}

This is another test case \cite{huang2011cures} which assesses the scheme for shock instability. This test case has complex flow features involving planar moving shock moving with incident Mach number, a diffracted shock wave around the corner and a strong expansion wave.  This strong expansion wave accelerates the flow and interacts with post shock fluid to further complicate the flow. Other distinct flow features are a slipstream and a contact surface. Numerical studies for this problem are given by Hillier \cite{hillier1991computation}. Godunov-type and Roe schemes are known to fail for this case \cite{quirk1994a}, as they produce oscillations at the planar shock.  As can be seen from the Fig. (\ref{FigShkDiffr}), the present method does not encounter any shock instability.  It can also be noted that slipstream is seen with first order simulation itself, whereas second order results also resolve contact surfaces.

\section{Conclusion}
A entropy-stable, hybrid numerical method for inviscid Euler equations is presented. This algorithm being a central solver, is free of complicated Riemann solvers and is not heavilty dependent on eigen-structure unlike upwind methods. This scheme uses the R-H condition-based numerical diffusion at the discontinuities and entropy conservation equation-based numerical diffusion in the smooth regions of the flow to capture grid aligned steady contact discontinuities and shock waves exactly, which is demonstrated through various benchmark problems. The idea of using numerical diffusion based on entropy conservation equation yields an entropy stable scheme for the smooth regions which is demonstrated using benchmark problems.  This hybrid scheme is also tested on strong shock test cases to prove its robustness and its accuracy is comparable to the best of the approximate Riemann solvers.  This entropy stable new algorithm is not prone to shock instabilities such as carbuncle shocks and avoids expansion shocks.  This new algorithm, MOVERS-LE, captures steady contact discontinuities and shock waves  exactly, which is one of the desired features for a scheme in simulating viscous flows, an aspect which is currently being pursued.

\bibliography{cfdreferences}

\end{document}